\documentstyle{article}

\newtheorem{theorem}{Theorem}[section]

\newtheorem{lemma}[theorem]{Lemma}
\newtheorem{definition}[theorem]{Definition}
\newtheorem{example}[theorem]{Example}

\newcommand{\defin}[1]
  {\begin{definition} {\rm #1} \end{definition}}
\newcommand{\examp}[1]
  {\begin{example} {\rm #1} \end{example}}

\def\QED{\quad\blackslug\lower 8.5pt\null}

  \newcommand{\htimes}{%
    \mathbin{\mathsurround0pt \mathchoice
      {\hbox{\vrule\negthinspace$\times$}}%
      {\hbox{\vrule\negthinspace$\times$}}%
      {\hbox{\vrule\negthinspace$\scriptstyle\times$}}%
      {\hbox{\vrule\negthinspace$\scriptscriptstyle\times$}}%
    }}

\begin{document}

{\Large \bf  Conformal and Grassmann structures}

\vspace*{3mm}

{\large M. A. Akivis}  

{\footnotesize \it             
Department of Mathematics,   Ben-Gurion University of the Negev,  
P.O. Box 653,  Beer Sheva 84105, 

Israel;  E-mail address: akivis@avoda.ict.ac.il       
}

\vspace*{3mm}

{\large   V.V. Goldberg} 

{\footnotesize \it             
Department of Mathematics,   New Jersey Institute of Technology,   University Heights, Newark, NJ 

07102, U.S.A.;  E-mail address: vlgold@numerics.njit.edu}

\vspace*{3mm}

{\footnotesize 
{\it Abstract:} 
The main results on the theory of conformal and almost Grassmann 
structures are presented. 
The common properties of  these structures and also the 
differences between them are outlined. In 
particular,  the structure groups of these structures and 
their differential prolongations are found.  A 
complete system of geometric objects of the almost Grassmann 
structure totally defining its geometric 
structure is determined. The vanishing of these objects 
determines a locally Grassmann manifold. It 
is proved that the integrable almost Grassmann structures are 
locally  Grassmann. The criteria of 
semiintegrability of almost Grassmann  structures is proved  in 
invariant form.

\vspace*{2mm}


\vspace*{2mm}

{\it MS classification}: 53A30, 53A40. 
}

\setcounter{equation}{0}

\setcounter{section}{-1}

\section{Introduction}  
The theory of conformal structures arose  in studying  
those properties of Riemannian and pseudo-Riemannian manifolds 
that remain invariant under conformal transformations of the 
metric. This theory was studied by many authors 
(see, for example, the paper [W 18] by Weyl, who first 
defined the tensor of conformal curvature of a Riemannian 
manifold, and the paper [C 23] by  Cartan, who introduced an 
$n$-dimensional space with a conformal connection).

Almost Grassmann manifolds were introduced by Hangan [H 66] 
as a generalization of the Grassmannian $G(m, n)$.  Hangan 
[H 66, 68] and Ishihara [I 72] studied mostly some special almost 
Grassmann manifolds, especially locally Grassmann manifolds.   
Later the almost Grassmann manifolds were studied 
by Goldberg [G 75], Mikhailov [M 78] and 
Akivis [A 82a] in connection with the construction of 
the theory of multidimensional webs. 
Goncharov [Go 87] considered the almost Grassmann manifolds  as 
generalized conformal structures. 

~Baston [B 91] constructed a theory of a general class of 
structures, called almost Hermitian symmetric (AHS) structures, 
which include conformal, projective, almost Grassmann, and 
quaternionic structures and for which the construction of the 
Cartan normal connection is possible. He constructed a tensor 
invariant for them and proved that its vanishing is equivalent to 
the structure being locally that of a Hermitian symmetric space.  
In [Go 87], the AHS structures have been studied from 
the point of view of cone structures. Bailey and 
Eastwood [BE 91] extended the theory of local twistors, which was 
known for four-dimensional conformal structures, to the  
almost Grassmann structures (they  called them the paraconformal 
structures).   
Dhooghe [Dh 94] considered the almost Grassmann structures (he 
called them Grassmannian structures) as subbundles of the second 
order frame bundle and constructed a canonical normal connection 
for these structures. The structure equations derived in [Dh 94] 
are very close to the structure equations of the almost Grassmann 
structures considered in the present paper.

Note that in the current paper we consider  the real theory 
of conformal and almost Grassmann structures while in [Go 87], 
[B 91] and [BE 91] their complex theory was studied.

Although some of the authors who studied almost Grassmann  
structures proved that an almost Grassmann  structure is 
a $G$-structure of finite type two (see [H 80] and [M 78]), 
no one of the authors went further than the 
development of the  first structure tensor. 

In the current paper we present the main results on the theory of 
conformal and almost Grassmann structures. The presentation is of 
survey nature, and as a rule, it does not contain complete proofs 
(most of them can be found in the papers [A 82a, b; 83, 85], 
[AK 93], [G 75], [H 80], [M 78], and in our book 
[AG 96]). However, the paper contains some new results, mostly 
related to the almost Grassmann  structures. In particular, 
in a third-order neighborhood we 
construct a complete system of geometric objects of the almost 
Grassmann  structure totally defining its geometric structure. 
The vanishing of these objects determines a locally Grassmann 
manifold. As for conformal structures, the complete  object of 
the almost Grassmann  structure is determined in a fourth order 
differential neighborhood.

In the theory of almost Grassmann  structures, integrable and 
semiintegrable structures play an important role. The integrable 
almost Grassmann structures are locally  Grassmann. The condition 
of semiintegrability of almost Grassmann  structures was first 
found in [M 78]. However, the proof in [M 78] has been done in a 
certain reduced frame bundle. Unlike the proof in [M 78], our 
proof in Section {\bf 4} is given in invariant form.

In the paper, we consider simultaneously the proper conformal 
structure $CO (n)$ and the  pseudoconformal structures 
$CO (p, q)$ of signature $(p, q)$. 

We find the common properties of conformal and almost Grassmann  
structures and also the differences between them. In particular, 
we find the structure groups of these structures and their 
differential prolongations. The structure group $G$ of the 
conformal structure $CO (p, q)$ is represented in the form 
${\bf SO} (p, q) \times {\bf H}$, and the structure group of the  
almost Grassmann structure in the form  
${\bf SL} (p) \times {\bf SL} (q) \times {\bf H}$, 
where  ${\bf SO} (p, q)$ is the special orthogonal group  of 
signature $(p, q)$;  ${\bf SL} (p), {\bf SL} (q)$ are the special 
linear groups of order $p$ and $q$, respectively; and ${\bf H}$ 
is the group of homotheties. For both structures, the prolonged 
group $G'$  is isomorphic to the semidirect product 
$G \htimes {\bf T} (n)$, where ${\bf T} (n)$ is the 
$n$-dimensional group  of parallel translations acting on $M$, 
and the $\htimes$ is the symbol of the semidirect product, but 
$n = p + q$ for the conformal structure and $n = pq$ for the 
almost Grassmann structure. 

The almost Grassmann  structure defines on the manifold $M$ two 
fiber bundles $E_\alpha$ and $E_\beta$. For the conformal 
structure, these fiber bundles arise only if $p = q = 2$. For the 
general almost Grassmann structure and the four-dimensional 
conformal structure, the first nonvanishing structure tensor 
splits into two subtensors, which are the structure tensors of 
fiber bundles $E_\alpha$ and $E_\beta$. 
For four-dimensional conformal and almost 
Grassmann structures the vanishing of each of 
these subtensors leads to integrability of the corresponding 
fiber bundle.

\section{Conformal Structures}
 
{\bf 1.} It is known that at any point $x$ of a pseudo-Euclidean 
space $R^n_q$ of signature $(p, q), \;\; p+q=n$, there is an 
isotropic cone $C_x (p, q)$  of second order with vertex at the 
point $x$. It  is also known that by compactification of the  
space $R^n_q$ one can construct a pseudoconformal space $C^n_q$ 
of the same signature. The compactification mentioned above is 
the enlargement of the space $R^n_q$ by the point at infinity, 
$y=\infty$,  and by the isotropic cone $C_y$ with  vertex 
at this point:
$$
C^n_q = R^n_q \bigcup \{C_y\}
$$

The space $C^n_q $ is a homogeneous space with the fundamental 
group 
$$
\renewcommand{\arraystretch}{1.3}
{\bf PO} (n+2, q+1) \cong \left\{ 
\begin{array}{ll}
{\bf SO} (n+2, q+1) \;\; \mbox{if} \;\; n \;\; \mbox{{\rm is odd}},\\
 {\bf O} (n+2, q+1)/{\bf Z}_2 \;\; \mbox{if} \;\; n \;\; 
\mbox{{\rm is even}}.
\end{array}
\right.
\renewcommand{\arraystretch}{1}
$$

Applying Darboux mapping, we can realize a pseudoconformal space 
$C^n_q$ on a hyperquadric $Q^n_q$ in a projective 
space $P^{n+1}$. After reduction to a sum of squares, the 
left-hand side of an equation of this hyperquadric will have 
$p+1$ positive and $q+1$ negative squares. Under the Darboux 
mapping, to the isotropic cones of the space  
$C^n_q$ there correspond 
the asymptotic cones of the hyperquadric $Q^n_q$ which 
are the intersections of the hyperquadric $Q^n_q$ with its 
tangent subspaces:
$$
C_x (p, q) = Q^n_q \cap T_x (Q^n_q), \;\; x \in Q^n_q.
$$
Note that for $q = 0$, the cone $C_x$ is imaginary, and for 
its consideration one must complexify the tangent space $C_x$, 
i.e. to enlarge it to the space ${\bf C} T_x = T_x \otimes 
{\bf C}$. 

Consider a family of moving frames in the space $C^n_q$, each of 
which is made up of  two points $A_0 = x$ and $A_{n+1}$ 
of general position and $n$ 
independent hyperspheres $A_i, \; i = 1, \ldots, n$, passing 
through these points. Under the Darboux mapping, this frame will 
pass into a point frame of a projective space $P^{n+1}$ such that 
its points $A_0$ and $A_{n+1}$ lie on the hyperquadric $Q^n_q$, 
and the points $A_i$, form a basis of the $(n-1)$-dimensional 
plane of intersection of two $n$-dimensional planes which are 
tangent to $Q^n_q$ at the points $A_0$ and $A_{n+1}$. If we 
denote by $(\;\;, \;\;)$ the scalar product of points in the 
space $P^{n+1}$, then the frame elements of the frames we have 
chosen satisfy the following analytical conditions: 
\begin{equation}\label{eq:1.1}
 (A_0, A_0) = 0, \;\; (A_{n+1}, A_{n+1}) = 0, \;\;
 (A_0, A_i) = 0, \;\; ( A_{n+1}, A_i) =0.  
\end{equation}
In addition, we normalize the points $A_0$ and $A_{n+1}$ by the 
condition 
\begin{equation}\label{eq:1.2}
 (A_0, A_{n+1}) = - 1.
\end{equation}
We will not assume that the hyperspheres $A_i$ are orthogonal 
and will write their scalar products in the form
\begin{equation}\label{eq:1.3}
 (A_i, A_j) =  g_{ij}, \;\; i, j = 1, \ldots , n.  
\end{equation}
(cf. [C 23] and [Y 39]). 
In the space $P^{n+1}$,  the hyperquadric $Q^n_q$ has the 
following equation with respect to the chosen frame: 
\begin{equation}\label{eq:1.4}
   g_{ij} x^i x^j - 2 x^0 x^{n+1} = 0.
\end{equation}
The quadratic form $g_{ij} x^i x^j$ is of signature $(p, q)$, 
i.e. its canonical form contains $p$ positive and $q$ negative 
squares.

We will write the equations of infinitesimal displacement of the 
moving frames  in the form 
\begin{equation}\label{eq:1.5}
dA_u = \omega_u^v A_v, \;\; u, v = 0, 1, \ldots , n+1,
\end{equation}
where $ \omega_u^v$ are differential 1-forms satisfying the 
structure equations of the space $P^{n+1}$:
\begin{equation}\label{eq:1.6}
d\omega_u^v = \omega_u^w \wedge \omega_w^v,
\end{equation}
which are the integrability conditions of equations (1.5). If we 
differentiate conditions (1.1)--(1.3), we obtain the following 
equations which the forms $\omega_u^v$ satisfy:
\begin{equation}\label{eq:1.7}
\renewcommand{\arraystretch}{1.5}
\begin{array}{ll}
\omega_0^{n+1} = \omega_{n+1}^0 = 0,   \;\;
\omega_0^0+ \omega_{n+1}^{n+1}  = 0, \\
\omega_i^{n+1} - g_{ij} \omega_0^j = 0,  \;\;
\omega_i^0 - g_{ij} \omega_{n+1}^j= 0 
\end{array}
\renewcommand{\arraystretch}{1}
\end{equation}
and 
\begin{equation}\label{eq:1.8}
dg_{ij} = g_{ik} \omega_j^k + g_{kj} \omega_i^k. 
\end{equation}

By equations (1.7), only the forms $\omega_0^0, \;\omega_0^i, \;
\omega_i^0$ and $\omega_i^j$ are linearly independent. In 
addition, the forms $\omega_i^j$ are connected with the 
fundamental tensor $g_{ij}$ by equations (1.8). 

The family of frames we have constructed in the space $C^n_q$ 
depends on $(n+1)^2$ parameters but this number does not 
coincide with the number $r$ of parameters on which the group 
of conformal transformations of the space $C^n_q$ depends. 
The latter number is equal to $(n + 1)^2 - \frac{1}{2} n (n + 1)$ 
where the subtrahend is equal to the number of independent among 
equations (1.8), i.e. $r = \frac{1}{2} (n + 1) (n + 2)$.  

Thus, the structure equations of the conformal  space $C^n_q$ 
can be written as 
\begin{equation}\label{eq:1.9}
\renewcommand{\arraystretch}{1.3}
\left\{
\begin{array}{ll}
\nabla g = 0, \\
d \omega = \kappa \wedge \omega - \theta \wedge \omega,\\
d \kappa = - \varphi \wedge \omega,\\
d \theta = \varphi \wedge \omega - \theta \wedge \theta + 
            (g \omega) \wedge (\varphi g^{-1}),\\
d \varphi = \varphi \wedge \kappa - \varphi \wedge \theta,
\end{array}
\right.
\renewcommand{\arraystretch}{1}
\end{equation}
where $g = (g_{ij}), \nabla g = (dg_{ij} - g_{ik} \omega_j^k  
- g_{kj} \omega_i^k), \omega = (\omega^i)$ (where $\omega^i = 
\omega_0^i$), 
$\kappa = \omega_0^0, \theta = (\omega_j^i), 
\varphi = (\omega_i^0)$, 
$d$ is the operator of exterior differentiation, and 
$\wedge$ is the symbol of exterior matrix multiplication. Note 
that in all exterior products of 1-forms occurring in equations 
(1.9)  multiplication is performed row by  column: for example, a 
detailed writing of second equation (1.9) has the following form:
$$
d \omega^i = \omega_0^0 \wedge \omega^i  
- \omega^i_j \wedge \omega^j. 
$$

{\bf 2.} A generalization of a conformal space is a conformal 
structure $CO (p, q)$ defined on a differentiable manifold $M$ of 
dimension $n$. 

Let $M$ be a differentiable manifold of dimension $n = p + q$ 
defined over reals ${\bf R}$. Consider the frame bundle 
consisting of vectorial frames $\{e_i\}$ and the conjugate 
coframe bundle $\{\omega^i\}$ consisting of 1-forms on $M$ such 
that 
$$
\omega^i (e_j) = \delta^i_j.
$$

Let $g$ be a nondegenerate quadratic form of signature $(p, q)$ 
(where $p+q = n$): 
\begin{equation}\label{eq:1.10}
g = g_{ij} \omega^i \omega^j, \;\;i, j, k = 1, \ldots, n. 
\end{equation}
The form $g$ defines  a {\em Riemannian metric} on $M$, and the 
tensor $g_{ij} = g_{ji}$ is its {\em metric tensor}. The 1-form 
$\omega = \{\omega^i\}$ is a vectorial form with its values in 
$T_x (M)$. This 1-form is defined in the first order frame bundle 
over $M$.  

A pair $(M, g)$ is called a {\em Riemannian manifold}. 

Two  Riemannian metrics $g$ and $\overline{g}$ given on the 
manifold $M$  are {\em conformally equivalent} if there exists a 
function $ \sigma (x) \neq 0, \; x \in M$, such that 
$$
 \overline{g} = \sigma g. 
$$

\defin{\label{def:1.1}
A {\em conformal structure} on the manifold  $M$ is the 
collection of all conformally equivalent Riemannian metrics given 
on $M$. 
}

We will denote such a structure by  $CO (p, q)$. Note that on the 
conformal structure $CO (p, q)$, the form (1.10) is relatively 
invariant, and on the Riemannian manifold $(M, g)$, it is 
absolutely invariant.

Note also that if $ \sigma (x) > 0$, then the forms $g$ and $\overline{g}$ have the same signature, and if $ \sigma (x) < 0$,   
then the form $\overline{g}$ is of signature $(q, p)$. Thus, 
{\em the structures $CO (p, q)$ and $CO (q, p)$ are equivalent}: 
$CO (p, q) \sim CO (q, p)$.

The equation $g = 0$ defines in ${\bf C} T_x (M) = T_x (M) \otimes {\bf C}$ a cone $C_x$ of second 
order and of signature $(p, q)$ which is called the 
{\em isotropic cone}: 
$$
C_x = \{\xi \in  {\bf C} T_x (M) | g (\xi, \xi) = 0\}.
$$

Conversely, it is easy to see that the structure {\em $CO (p, q)$ 
is defined on a real manifold $M$ by a fibration 
of cones $C_x$ of second order and 
of signature $(p, q)$}. Thus, 
$$
 CO (p, q) = \Bigl(M, \displaystyle \bigcup_{x \in M} 
C_x (p, q) | C_x \subset {\bf C} T_x (M)\Bigr).
$$

The structure group $G$ of the structure  $CO (p, q)$ is locally 
isomorphic to the subgroup of ${\bf GL} (n, {\bf R})$ leaving 
invariant the cone $C_x$:
$$
 G \cong {\bf SO} (p, q)  \times {\bf H}, \;\; p+q = n, 
$$
where ${\bf SO} (p, q)$ is the special $n$-dimensional 
pseudoorthogonal group of signature $(p, q)$ (the connected 
component of the unity of the pseudoorthogonal group 
${\bf O} (p, q)$),  ${\bf H}$ is the group of homotheties, and 
$\cong$ is the symbol of local isomorphism.  

In ${\bf C} T_x (M)$, the action of the group $G$ is as follows:
$$
\begin{array}{ll}
G = \{\gamma \in {\bf SO} (p, q) \times {\bf H}, \xi \in 
{\bf C} T_x (M) |\\ 
 \gamma T_x = T_x, \gamma (\overline{\xi}) 
= \overline{\gamma (\xi)}, \gamma C_x = C_x\}.
\end{array}
$$

It follows that {\em the structure $CO (p, q)$ is a $G$-structure 
of first order on  $M$ with  the structure group $G$}.  

{\bf 3.} Let us consider some particular cases and examples of 
$CO (p, q)$-structures. 

First, if  $p = n$ and $q = 0$, then we have the {\em proper 
conformal structure} $CO (n, 0) = CO (n)$. In this case, the cone 
$C_x$ is imaginary,  and structural group 
$G \cong {\bf O} (n) \times {\bf H}$. If  $0 < q < n$, then we 
have a  pseudoconformal structure. For this structure, the cone 
$C_x$ is real. In particular, if $p = 1$ and $q = n - 1$, then we 
obtain a pseudoconformal structure of {\em Lorentzian type}, and 
if $2 \leq p \leq n - 2$, the $CO (p, q)$-structure is called 
{\em ultrahyperbolic}.

\examp{\label{examp:1.2} 
Since there is  an isotropic cone $C_x$ of signature $(p, q)$ 
at any point $x$ of  
the pseudoconformal space $C^n_q $ defined in subsection {\bf 1},  
 this space carries a $CO (p, q)$-structure. 
We will call the $CO (p, q)$-structure associated with the space 
$C^n_q$ {\em conformally flat}.

It is obvious that the quadric $Q^n_q$, which arises 
if one apply the Darboux mapping to the space $C^n_q$, carries 
also a conformally flat $CO (p, q)$-structure defined by the 
asymptotic cones of $Q^n_q$ ($C_x = Q^n_q \bigcap T_x (Q^n_q)$) (see subsection {\bf 1}).
}

\examp{\label{examp:1.3} Suppose that $n = 4$. Then there are 
three $CO (p, q)$-structures: $CO (4, 0)$-  (or $CO (4)$-), 
$CO (1, 3)$- and $CO (2, 2)$-structures. 

On the $CO (2, 2)$-structure, the cones $C_x$ carry two families 
of real plane generators, $\alpha$- and $\beta$-planes, that form 
the {\em isotropic distributions} $E_\alpha$ and $E_\beta$ on 
$M$. Thus, 
$G \cong {\bf SL} (2) \times {\bf SL} (2) \times {\bf H}$. 

On the $CO (1, 3)$-structure, the cones $C_x$ carry a family of 
real rectilinear generators and two families of complex plane 
generators;  $E_\beta = \overline{E}_\alpha$ and 
$G \cong {\bf SL} (2, {\bf C}) \times {\bf H} 
\cong \overline{{\bf SL} (2, {\bf C})} \times {\bf H}$,  and 
the groups ${\bf SL} (2, {\bf C})$ and $ \overline{{\bf SL}(2, {\bf C})}$ act concordantly on $E_\alpha$ and $E_\beta$.  
Note that in general relativity, $C_x$ are light cones.

On the $CO (4)$-structure, the cones $C_x$ are pure imaginary, 
$E_\alpha$ and $E_\beta$ are self-conjugate: 
$\overline{E}_\alpha = E_\alpha; \; 
\overline{E}_\beta = E_\beta$; 
 $G \cong {\bf SU} (2) \times {\bf SU} (2) \times {\bf H}$, 
where ${\bf SU} (2)$ is the special two-dimensional unitary 
group.
}

{\bf 4.}
The 1-form $\omega = \{\omega^i\}$ is defined in the bundle  
${\cal R}^1 (M)$ of frames of first order. In addition to this 
form, one can invariantly define a matrix 1-form  
$\theta = (\omega_j^i)$  and a scalar form 
$\kappa$ in the  bundle  ${\cal R}^2 (M)$ of frames of second 
order, and a covector 1-form  $\varphi = (\omega_i)$ in the  
bundle  ${\cal R}^3 (M)$ of frames of third order (see [AK 93]). 
All these forms satisfy the following structure equations: 
\begin{equation}\label{eq:1.11}
\renewcommand{\arraystretch}{1.3}
\left\{
\begin{array}{ll}
\nabla g = 0,\\
d \omega = \kappa \wedge \omega 
- \theta \wedge \omega, \\
 d \kappa = - \varphi  \wedge \omega,\\
d \theta = \varphi  \wedge \omega - \theta 
\wedge \theta + (g \omega) \wedge 
(\varphi g^{-1}) + \Theta, \\
d \varphi = \varphi  \wedge \kappa - \varphi \wedge \theta 
+ \Phi, 
\end{array}
\right.
\renewcommand{\arraystretch}{1}
\end{equation}
where, as in formulas (1.9), 
$\nabla g = (dg_{ij} - g_{ik} \omega_j^k 
- g_{kj} \omega_i^k)$, and the 2-forms $\Theta = (\Theta_j^i)$ 
and $\Phi = (\Phi_i)$ are the curvature 2-forms. Note that 
for $\Theta = \Phi = 0$, the  structure equations (1.11) coincide 
with the structure equations (1.9) of the pseudoconformal space 
$C^n_q$. Thus, {\em the space $C^n_q$ 
carries the curvature-free 
$CO (p, q)$-structure, i.e. the space $C^n_q$ 
is conformally flat.}

We will now clarify the geometric meaning for the 1-forms 
$\theta, \kappa$, and $\varphi$. 

First, we note that the first of  equations (1.11) is the 
condition for the cone $C_x$ to be invariant. As to other 
equations of (1.11), if $\omega^i = 0$, then they become 
the structure equations of a 
pseudoconformal space $C^n_q$ in which also $\omega^i = 0$ being 
set (i.e. a point $A_0$ is fixed):
\begin{equation}\label{eq:1.12}
 d \kappa = 0, \;\;
d \theta  = - \theta \wedge \theta, \;\; 
d \varphi = \varphi  \wedge \kappa - \varphi \wedge \theta. 
\end{equation}

 From equations (1.12) it follows that 

\begin{description}

\item[(i)] The form $\kappa$ is an invariant form of the group 
${\bf H}$ of homotheties acting in the space $T_x (M)$. 

\item[(ii)]  The form $\theta$ is an invariant forms of the 
special pseudoorthogonal group ${\bf SO} (p, q)$ leaving the 
points $x$ and $y$ invariant.

\item[(iii)] The forms $\kappa$ and $\theta$ are invariant forms 
of the structural group $G$ of the pseudoconformal structure 
$CO (p, q)$ whose transformations leave cone $C_x (p, q)$ 
invariant. The  group $G$ is isomorphic to the direct product 
${\bf SO} (p, q) \times {\bf H}$: 
\begin{equation}\label{eq:1.13}
G \cong {\bf SO} (p, q) \times {\bf H}.
\end{equation}

\item[(iv)] 
The 1-forms $\theta, \kappa$, and $\varphi$ are invariant forms 
of the group $G'$ which is a differential prolongation of $G$. 
The group $G'$ is the group of motions of the space $(C^n_q)_x$ 
which is the compactification of the tangent space $T_x (M)$: 
$(C^n_q)_x = T_x (M) \bigcup C_y $, and  $(C^n_q)_x$ is referred  to a frame $\{x, y, a_i\}$ where $a_i$ are hyperspheres 
passing through the points $x$ and $y$.

\item[(v)] The covector 
$\varphi$  is  an invariant form of  the group of translations 
 ${\bf T} (p+q)$ of the pseudo-Euclidean space $R^n_q 
= C^n_q \backslash C_x$ whose 
transformations move the point  $y$. 
\end{description}

 Thus, the group $G'$ is isomorphic to the semidirect 
product $G \htimes {\bf T} (p+q)$, i.e. 

\begin{equation}\label{eq:1.14}
G' \cong G \htimes {\bf T} (p+q) \cong
 ({\bf SO} (p, q) \times {\bf H}) \htimes {\bf T} (p+q), 
\end{equation}
and is the group of motions of the space $R^n_q$.

Since the group $G'$ does not admit further prolongations,  
{\em a pseudoconformal structure $CO (p, q)$ is a $G$-structure 
of finite type two} (see [S 64] for definition of finite type). 

The structure equations of the $CO (p, q)$-structure in the form 
(1.11) can be found in [C 23], [Ga 89], [A 85] and [AK 93].  
Cartan [C 23]  called equations (1.11) the structure equations of 
the normal conformal connection associated with the quadratic 
form $g$.  Note also that while in [C 23] only proper conformal structures were considered, in  [C 22a, b] the 
pseudoconformal structures were studied as well.

{\bf 5.} The forms 
\begin{equation}\label{eq:1.15}
\Theta^i_j = b^i_{jkl} \omega^k \wedge \omega^l, \;\;
\Phi_i = c_{ijk} \omega^j \wedge \omega^k
\end{equation}
appearing in the last two equations of (1.11) are called the 
{\em curvature forms} of the $CO (p, q)$-structure. The 
quantities $b^i_{jkl}$ are the components of the  {\em tensor 
of conformal curvature} (the Weyl tensor) $b = \{b^i_{jkl}\}$ of 
the $CO(p, q)$-structure, and the quantities $c_{ijk}$  together 
with $b^i_{jkl}$  constitute a homogeneous geometric object 
$(b, c) = \{b^i_{jkl}, c_{ijk}\}$. The tensor  $b^i_{jkl}$ is 
defined in a third-order differential neighborhood of a point 
$x \in M$. It satisfies all conditions which the curvature tensor 
of a Riemannian manifold satisfies, and in addition, the tensor 
$b^i_{jkl}$ is trace-free:
\begin{equation}\label{eq:1.16}
b^i_{jki} = 0
\end{equation}
(see, for example, the book [E 26]).

The geometric object $(b, c)$ is defined in a fourth differential 
neighborhood of a point $x \in M$. If $n \geq 4$, the quantities 
$c_{ijk}$ are expressed linearly in terms of the covariant 
derivatives of the tensor $b^i_{jkl}$:
\begin{equation}\label{eq:1.17}
c_{ijk} = - \frac{1}{n-3} b^m_{ijk,m}.
\end{equation}

Thus, if  $n \geq 4$ and $b = 0$,  then $c = 0$, 
and the $CO (p, q)$-structure is locally flat. Note 
that if $n = 3$, then  $b \equiv 0$, 
and the condition $c = 0$ is the 
condition for the $CO (p, q)$-structure to be locally flat. 

Since the second equation of (1.11) does not have an exterior 
quadratic form of type (1.15), a {\em $CO (p, q)$-structure is 
torsion-free}.

It is proved in [A 83] (see also [AK 93]) that for the 
$CO (2, 2)$-structure, the tensor $b$ of conformal curvature 
splits into two subtensors $b_\alpha$ and   $b_\beta$: 
$b = b_\alpha \dot{+} b_\beta$, and $b_\alpha$ and $b_\beta$ are 
the curvature tensors of 
the fibre bundles $E_\alpha$ and $E_\beta$ (see subsection 
{\bf 3}). Each of the subtensors $b_\alpha$ and $b_\beta$ has 
five independent components. 

 Since for the 
$CO (1, 3)$-structure, the fiber bundles $E_\alpha$ 
and $E_\beta$ are complex conjugates, its 
curvature tensor admits two complex conjugate representations 
$b_\alpha$ and $b_\beta$, $b_\beta = \overline{b}_\alpha$, 
which themselves are the curvature tensors of 
the fiber bundles $E_\alpha$ and $E_\beta$. For the 
$CO (4)$-structure, we again have the splitting 
$b = b_\alpha  \dot{+} b_\beta$, 
and the subtensors $b_\alpha$ and $b_\beta$ are self-conjugates:
$b_\alpha = \overline{b}_\alpha$ and 
$b_\beta = \overline{b}_\beta$.

For $p + q = 4$, a conformal $CO (p, q)$-structure 
 for which the tensor $b_\alpha$ or 
$b_\beta$ vanishes is called {\em conformally semiflat}. 
If both these 
tensors vanish, the structure is called {\em conformally flat}. 
The conditions which the tensors $b_\alpha$ and $b_\beta$ 
satisfy show that {\em only 
the $CO (2, 2)$- and $CO (4)$-structures can be  conformally 
semiflat}, and that {\em the $CO (1, 3)$-structure cannot be 
conformally semiflat without being conformally flat}.

It is easy to see that {\em the  conformally flat  
$CO (2, 2)$-structure  is locally isomorphic to 
the structure of the pseudoconformal space $C^4_2$}.

Finally note that for the $CO (4)$-structure the notion of local 
semiflatness is connected with the notions of self-duality 
and anti-self-duality introduced in [AHS 78].

\section{Grassmann Structures}

\setcounter{equation}{0}
 
{\bf 1.} 
Let $P^n$ be an $n$-dimensional projective space.  The set of 
$m$-dimensional subspaces $P^m \subset P^n$ is called the {\em 
Grassmann manifold}, or the {\em Grassmannian}, and is denoted by 
the symbol $G (m, n)$. It is well-known that the Grassmannian is 
a differentiable manifold, and that its dimension is equal to 
$\rho = (m + 1)(n - m)$. It will be convenient for us to set 
$p = m + 1$ and $q = n - m$. Then we have $n = p + q - 1$.
 
Let a subspace $P^m = x$ be an element of the Grassmannian 
$G (m, n)$. With any subspace  $x$, we associate a family  of 
projective point frames $\{A_u\}, u = 0, 1, \ldots n$, such that 
the vertices $A_\alpha, \alpha = 0, 1, \ldots , m$, of its frames 
lie in the plane $P^m$, and the points $A_i, i = m+1, \ldots, n$, 
lie outside $P^m$ and together with the points $A_\alpha$ make up 
the frame $\{A_u\}$ of the space  $P^n$. 

We will write the equations of infinitesimal displacement of the 
moving frames we have chosen in the form:
\begin{equation}\label{eq:2.1}
d A_u = \theta_u^v A_v, \;\; u, v = 0, \ldots, n.
\end{equation}
Since the fundamental group  of the space $P^n$ is locally 
isomorphic to the group ${\bf SL} (n+1)$, the forms 
$\theta_u^v$ are connected by the relation 
\begin{equation}\label{eq:2.2}
\theta_u^u = 0. 
\end{equation}
The structure equations of the space $P^n$ have the form
\begin{equation}\label{eq:2.3}
d \theta_u^v = \theta_u^w \wedge \theta_w^v. 
\end{equation}
By (2.3), the exterior differential of the left-hand side of 
equations (2.1) is identically equal to 0, and hence the system 
of equations (2.1) is completely integrable.

By (2.1), we have 
$$
d A_\alpha = \theta_\alpha^\beta A_\beta + \theta_\alpha^i A_i.
$$
It follows that the 1-forms $\theta_\alpha^i$ are basis forms of 
the Grassmannian. These forms are linearly independent, and their 
number is equal to $\rho = (m+1)(n-m) = p q$, i.e. it 
equals the dimension of the 
Grassmannian $G (m, n)$. We will assume  that the integers $p$ 
and $q$ satisfy the inequalities $p \geq 2$ and $q \geq 2$, since 
for $p = 1$, we have $m = 0$, and the Grassmannian $G (0, n)$ is 
the projective space $P^n$, and for $q=1$, we have $m = n - 1$, 
and the Grassmannian $G (n-1, n)$ is isomorphic to the dual 
projective space $(P^n)^*$.

Let us rename the basis forms by setting $\theta_\alpha^i = 
\omega_\alpha^i$ and find their exterior differentials:
\begin{equation}\label{eq:2.4}
d \omega_\alpha^i = \theta_\alpha^\beta \wedge \omega_\beta^i 
+ \omega_\alpha^j \wedge \theta^i_j. 
\end{equation}
Define the trace-free forms 
\begin{equation}\label{eq:2.5}
\omega_\alpha^\beta = \theta_\alpha^\beta 
- \frac{1}{p} \delta_\alpha^\beta \theta^\gamma_\gamma, \;\;
\omega^i_j = \theta^i_j  
- \frac{1}{q} \delta_i^j  \theta^k_k,
\end{equation}
satisfying the conditions
\begin{equation}\label{eq:2.6}
\omega_\alpha^\alpha = 0, \;\;\; \omega_i^i = 0.
\end{equation}
 Eliminating the forms $\theta_\alpha^\beta $ and $\theta_i^j$ 
from equations (2.4), we find that
\begin{equation}\label{eq:2.7}
d \omega_\alpha^i = \omega_\alpha^\beta \wedge \omega_\beta^i 
+ \omega_\alpha^j \wedge \omega^i_j + \kappa \wedge \omega_\alpha^i, 
\end{equation}
where $\kappa = \displaystyle \frac{1}{p} \theta^\gamma_\gamma 
- \displaystyle \frac{1}{q} \theta^k_k$, or by (2.2), 
\begin{equation}\label{eq:2.8}
\kappa = \Bigl(\displaystyle \frac{1}{p} + \displaystyle  
\frac{1}{q}\Bigr)   \theta^\gamma_\gamma.
\end{equation}

Setting
\begin{equation}\label{eq:2.9}
 \omega_i^\alpha = - \Bigl(\displaystyle \frac{1}{p} + \displaystyle  \frac{1}{q}\Bigr) \theta_i^\alpha    
\end{equation}
and taking the exterior derivatives of equations (2.7) and (2.8), 
we obtain 
\begin{equation}\label{eq:2.10}
\renewcommand{\arraystretch}{1.3}
\left\{
\begin{array}{ll}
d \omega_\alpha^\beta = \omega_\alpha^\gamma \wedge 
 \omega_\gamma^\beta + \displaystyle \frac{q}{p+q} 
\omega_\gamma^k \wedge 
(\delta_\alpha^\beta \omega_k^\gamma - p \delta_\alpha^\gamma 
\omega_k^\beta), \\
d \omega_j^i = \omega_j^k \wedge \omega_k^i + 
\displaystyle \frac{p}{p+q} (\delta^i_j \omega_k^\gamma 
- q \delta^i_k \omega_j^\gamma) 
\wedge \omega_\gamma^k,  
\end{array}
\right.
\renewcommand{\arraystretch}{1}
\end{equation}
and 
\begin{equation}\label{eq:2.11}
 d \kappa = \omega_i^\alpha \wedge \omega^i_\alpha.    
\end{equation}
Exterior differentiation of equations (2.9) gives 
\begin{equation}\label{eq:2.12}
 d \omega_i^\alpha  = \omega_i^j \wedge \omega_j^\alpha 
+ \omega_i^\beta \wedge \omega_\beta^\alpha + 
\omega_i^\alpha \wedge \kappa.    
\end{equation}
Finally,  exterior differentiation of equations (2.12) leads to 
identities.

Thus, the structure equations of the Grassmannian $G (m, n)$ take 
the form (2.7), (2.10),  (2.11), and (2.12). This system of 
differential equations is {\em closed} in the sense that its 
further exterior differentiation leads to identities. 

If we fix a subspace $x = P^m \subset P^n$, then we obtain 
$\omega^i_\alpha = 0$, and equations (2.10) and (2.11) become
\begin{equation}\label{eq:2.13}
d \pi_\alpha^\beta = \pi_\alpha^\gamma \wedge 
 \pi_\gamma^\beta, \;\; d \pi_j^i = \pi_j^k \wedge \pi_k^i, \;\;  
d \pi = 0,  
\end{equation}
where  $ \pi = \kappa (\delta), 
\pi_\alpha^\beta = \omega_\alpha^\beta (\delta), 
\pi_j^i = \omega_j^i (\delta)$, 
and $\delta$ is the operator of differentiation with respect to 
the fiber parameters of the second order frame bundle associated 
with the Grassmannian $G (m, n)$. Moreover, the forms 
$\pi_\alpha^\beta$ and $\pi_j^i$ satisfy  equations similar to 
equations (2.6), i.e. these forms are trace-free. The forms 
$\pi_\alpha^\beta$ are invariant forms of the group 
${\bf SL} (p)$ which is locally isomorphic to the group of 
projective transformations of the subspace $P^m$. The forms 
$\pi_j^i$ are invariant forms of the group ${\bf SL} (q)$ which 
is locally isomorphic to the group of projective transformations 
of the bundle  of $(m+1)$-dimensional subspaces of the space 
$P^n$ containing $P^m$. The form $\pi$ is an invariant form of 
the group ${\bf H} = {\bf R}^* \otimes \mbox{{\rm Id}}$ of the 
space $P^n$ with center at $P^m$; here ${\bf R}^*$ is the 
multiplicative group of real numbers. 

The direct product of these three groups is the structural group 
$G$ of the Grassmann manifold $G (m, n)$:
\begin{equation}\label{eq:2.14}
G = {\bf SL} (p) \times {\bf SL} (q) \times {\bf H}.
\end{equation}

Finally, the forms $\pi_i^\alpha = \omega_i^\alpha (\delta)$, 
which by (2.12) satisfy the structure equations 
\begin{equation}\label{eq:2.15}
d \pi_i^\alpha = \pi_i^j \wedge  \pi_j^\alpha + \pi_i^\beta 
\wedge \pi_\beta^\alpha + \pi_i^\alpha \wedge \pi,   
\end{equation}
are also fiber forms on the Grassmannian $G (m, n)$ but unlike 
the forms $\pi_\alpha^\beta,\; \pi_i^j$ and $\pi$, they are 
connected with the third-order frame bundle of the Grassmannian 
$G (m, n)$. 

The forms  $\pi_\alpha^\beta, \;\pi_j^i, \pi$, and 
$\pi_i^\alpha$, satisfying the structure equations 
(2.13) and (2.15), are invariant forms of the  group 
\begin{equation}\label{eq:2.16}
G' = G \htimes {\bf T} (pq)
\end{equation}
arising under the differential prolongation 
of the structure group $G$ of the Grassmannian $G (m, n)$. 
The group $G'$ is the group of motions of 
an $(n - m - 1)$-quasiaffine space 
$A^n_{n-m-1}$ (see [R 59]) which is 
a projective space $P^n$ with a fixed  $m$-dimensional subspace 
$P^m = A_0 \wedge A_1 \wedge \ldots \wedge A_m$ and 
the generating element 
$P^{n-m-1} = A_{m+1} \wedge \ldots \wedge A_n$. 
The dimension of the space 
$A^n_{n-m-1}$ coincides with the 
dimension  of the Grassmannian $G (n - m - 1, n)$, and this 
dimension is the same as the dimension of the Grassmannian 
$G (m, n)$: $\rho = (m + 1)(n - m)$. 
The forms  $\pi_i^\alpha$ are invariant forms of the  group  
${\bf T} (pq)$  of parallel translations of 
the space $A^n_{n-m-1}$, and  the group $G$  is the 
stationary subgroup  of its element $P^{n-m-1}$.

In the index-free notation, the structure equations (2.7) and 
(2.10)--(2.12) of the Grassmannian $G (m, n)$ can be written as 
follows:
\begin{equation}\label{eq:2.17}
\renewcommand{\arraystretch}{1.3}
\left\{
\begin{array}{ll}
 d\omega =  \kappa \wedge \omega -  \omega  \wedge \theta_\alpha 
-  \theta_\beta \wedge \omega, \\
d \theta_\alpha  + \theta_\alpha \wedge \theta_\alpha  
= \displaystyle \frac{q}{p+q} \Bigl[- I_\alpha 
\mbox{{\rm tr}}\; (\varphi \wedge \omega) 
+ p \varphi \wedge \omega\Bigr],  \\
d \theta_\beta  + \theta_\beta \wedge \theta_\beta  
= \displaystyle \frac{p}{p+q} \Bigl[- I_\beta 
\mbox{{\rm tr}}\; (\varphi \wedge \omega) 
+ q \omega \wedge \varphi\Bigr],  \\
d \kappa  =  \mbox{{\rm tr}} \; (\varphi \wedge \omega), \\
 d \varphi  + \theta_\alpha \wedge \varphi  + \varphi \wedge 
\theta_\beta  + \kappa \wedge  \varphi = - (a \varphi) 
\wedge \omega,
\end{array}
\right.
\renewcommand{\arraystretch}{1}
\end{equation}
where $\omega = (\omega_\alpha^i)$ is the matrix 1-form 
defined in the first-order fiber bundle;  
 $\kappa$ is a  scalar 1-form; 
$\theta_\alpha = (\omega_\beta^\alpha)$ and 
$\theta_\beta = (\omega_j^i)$ are the matrix 1-forms 
defined in a second order fiber bundle for which 
$$
\mbox{{\rm tr}} \; \theta_\alpha = 0, \;\; 
\mbox{{ \rm tr}} \; \theta_\beta = 0;
$$
$\varphi = (\omega_i^\alpha)$ is a matrix 1-form defined in a 
third-order fiber bundle; and
 $I_\alpha = (\delta_\alpha^\beta)$ and 
$I_\beta = (\delta_i^j)$ are the unit tensors of orders $p$ and 
$q$, respectively. 

Along with the Grassmannian $G (m, n)$, in the space $P^n$ 
one can consider the dual manifold $G (n - m - 1, n)$. Its base 
forms are the forms $\omega_i^\alpha$, and its geometry is identical to that of the Grassmannian $G (m, n)$.

{\bf 2.} With the help of Grassmann coordinates, the Grassmannian 
$G (m, n)$ can be mapped onto a smooth algebraic variety 
$\Omega (m, n)$ of dimension $\rho = pq$ embedded into a 
projective space $P^N$ of dimension 
$N =  {p+q \choose p} - 1$.

Suppose that $x = A_0 \wedge A_1 \wedge \ldots \wedge A_m $ 
is a point of the variety $\Omega (m, n)$. Then 
\begin{equation}\label{eq:2.18}
d x = \tau x + \omega_\alpha^i e^\alpha_i,
\end{equation}
where $\tau = \omega_0^0 + \ldots + \omega_m^m$ and 
$$
e_i^\alpha = A_0 \wedge \ldots  \wedge A_{\alpha-1} 
\wedge A_i  \wedge A_{\alpha+1}  \wedge \ldots  \wedge A_m, 
$$
and the points $e_i^\alpha$ together with the point $x$ determine 
a basis in the tangent subspace $T_x (\Omega)$. The second 
differential of the point $x$ satisfies the relation
\begin{equation}\label{eq:2.19}
d^2 x \equiv \sum_{\alpha<\beta, i< j} 
(\omega_{\alpha}^{i} \omega_{\beta}^{j} - \omega_{\alpha}^{j} \omega_{\beta}^{i})  e^{\alpha \beta}_{ij} \pmod{T_x (\Omega)},
\end{equation}
where  
$$
 e^{\alpha \beta}_{ij} =  A_0 \wedge \ldots  \wedge A_{\alpha-1} 
 \wedge A_i  \wedge A_{\alpha+1}  \wedge \ldots \wedge 
A_{\beta-1} \wedge A_j  \wedge A_{\beta+1}  \wedge \ldots \wedge 
A_m 
$$
are points of the space $P^N$ lying on the variety 
$\Omega (m, n)$. The quadratic forms 
\begin{equation}\label{eq:2.20}
\omega_{\alpha\beta}^{ij} = \omega_{\alpha}^{i} 
\omega_{\beta}^{j} - \omega_{\alpha}^{j} \omega_{\beta}^{i}  
\end{equation}
are the second fundamental forms of the variety $\Omega \subset 
P^N$. 

The equations $\omega_{\alpha\beta}^{ij} = 0$ determine the cone 
of asymptotic directions of the variety $\Omega$ at a point 
$x \in \Omega$. The equations of this cone can be written as 
follows:
\begin{equation}\label{eq:2.21}
{\rm rank}\;\; (\omega^i_\alpha) = 1.
\end{equation}
In view of (2.21), parametric equations of this cone have 
the form 
\begin{equation}\label{eq:2.22}
\omega_{\alpha}^{i} = t_{\alpha} s^{i}, \;\; 
\alpha = 0, 1, \ldots , m;\;\; i = m+1, \ldots , n. 
\end{equation}

If we consider a projectivization of this cone, then $t_\alpha$ 
and $s^i$ can be taken as homogeneous coordinates of projective 
spaces $P^m$ and $P^{n-m-1}$. Thus such a projectivization is an 
embedding of the direct product $P^{p-1} \times P^{q-1}$ into a 
projective space $P^{\rho - 1}$ of dimension ${\rho - 1}$. Such 
an embedding is called  the {\em Segre variety} and is denoted by $S (p-1, q-1)$. This is the reason that the cone of asymptotic 
directions of the variety $\Omega$ determined by equations (2.21) 
is called the {\em Segre cone}. This cone is denoted by 
$SC_x (p, q)$ since it carries two families of plane generators 
of dimensions $p$ and $q$.  Plane generators from different 
families of the cone $SC_x (p, q)$ have a common straight line. 
It is possible to prove that the cone $SC_x (p, q)$ is the 
intersection of the tangent subspace $T_x (\Omega)$ and the 
variety $\Omega$:
$$
C_x (p, q) = T_x (\Omega) \cap \Omega.
$$

The differential geometry of Grassmannians was studied in detail 
in the paper [A 82b].

\section{Almost Grassmann Structures}

\setcounter{equation}{0}
 
{\bf 1.} 
Now we can define the notion of an almost Grassmann structure. 

\defin{\label{def:3.1} Let $M$ be a differentiable manifold of 
dimension $p q$, and let $SC (p, q)$ be a differentiable 
fibration of Segre cones with the base $M$ such that 
$SC_x (M) \subset T_x (M), \;\; x \in M$. The pair 
$(M, SC (p, q))$ is said to be  an {\it almost Grassmann 
structure} and is denoted by $AG (p-1, p+q-1)$. The manifold 
$M$ endowed with  such a structure is said to be an {\it almost 
Grassmann manifold}.
}

As was the case for Grassmann structures, the almost Grassmann 
structure $AG (p-1, p+q-1)$ is equivalent to the structure 
$AG (q-1, p+q-1)$ since both of these structures are generated on 
the manifold $M$ by a differentiable family of  Segre cones 
$SC_x (p, q)$.

 Let us consider some examples. 

\examp{\label{examp:3.2}
The main example of an almost Grassmann structure is the almost 
Grassmann structure associated with the Grassmannian $G (m, n)$.  
As we saw, there is a field of Segre cones $SC_x (p, q) 
= T_x (\Omega) \bigcap \Omega, \; x \in \Omega$, where $p = m+1$ 
and $q = n-m$, which defines  an almost Grassmann structure.
}
 
\examp{\label{examp:3.3}
Consider a pseudoconformal $CO (2, 2)$-structure on a 
four-di\-men\-sional manifold $M$. The isotropic cones $C_x$ of 
this structure carry two families of plane generators. Hence, 
these cones are Segre cones $SC_x (2, 2)$. Therefore, a 
pseudoconformal $CO (2, 2)$-structure is an almost Grassmann 
structure $AG (1, 3)$. 

If we complexify the four-dimensional tangent subspace 
$T_x (M^4)$ and consider Segre cones with complex generators, 
then  conformal $CO (1, 3)$- and $CO (4, 0)$-structures can also 
be considered as complex almost Grassmann structures of the same 
type $AG (1, 3)$. However, in this paper, we will consider only 
real almost Grassmann structures.
}

Almost Grassmann structures arise also in the study of 
multidimensional webs (see  [AS 92] and [G 88]). 

\examp{\label{examp:3.4}
Consider a three-web formed on a manifold $M^{2q}$ of dimension 
$2q$ by three foliations $\lambda_u, \; u = 1, 2, 3$, of 
codimension $q$ which are in general position (see [A 69] 
and [AS 92]). 

Through any point $x \in M^{2q}$, there pass three leaves 
${\cal F}_u$ belonging to the foliations $\lambda_u$. In the 
tangent subspace $T_x (M^{2q})$, we consider three subspaces 
$T_x ({\cal F}_u)$ which are tangent to ${\cal F}_u$ at the point 
$x$. If we take the projectivization of this configuration with  
center at the point $x$, then we obtain a projective space 
$P^{2q-1}$ of dimension $2q-1$ containing three subspaces of 
dimension $q-1$ which are in general position. These three 
subspaces determine a Segre variety $S (1, q - 1)$, and the 
latter variety is the directrix for a Segre cone 
$SC_x (2, q) \subset T_x (M^{2q})$. Thus, on $M^{2q}$, a field of 
Segre cones  arises, and this field determines an almost 
Grassmann structure on $M^{2q}$. 

The structural group   of the web $W (3, 2, q)$ 
is smaller than that  of the induced  almost Grassmann structure, 
since transformations of this group must keep invariant the 
subspaces $T_x ({\cal F}_u)$. Thus, the structural group of 
the three-web is the group ${\bf GL} (q)$. 
}

\examp{\label{examp:3.5} 
 Consider a $(p+1)$-web $W(p+1, p, q) 
= (M; \lambda_1, \ldots, \lambda_{p+1})$ formed on a 
differentiable manifold $M$ of dimension $pq$ by $p+1$ foliations 
$\lambda_u, \; u=1, \ldots, p+1$, of dimension $q$ which are in 
general position on $M$ (see [G 73] or [G 88]).

As in Example {\bf 3.4}, the tangent spaces $T_x ({\cal F}_u)$ 
define the cone $SC_x (p, q) \supset T_x ({\cal F}_u)$, and 
the field of these cones defines an almost Grassmann structure 
$AG (p-1, p+q-1)$ on $M$.

The structural group of the web $W (p+1, p, q)$ is the same group 
$ G ={\bf GL} (q)$ as for the web $W (3, 2, q)$, and this group 
does depend on $p$.
}

{\bf 2.} 
The structural group of the almost Grassmann structure is a 
subgroup of the general linear group ${\bf GL}(pq)$ of 
transformations of the space $T_{x} (M)$, which leave the cone 
$SC_{x}(p, q) \subset T_x (M)$ invariant. We denote this group by 
$G = {\bf GL} (p, q)$. 

To clarify the structure of this group, in the tangent space  
$T_x (M)$, we consider a family of frames  
$\{e_{i}^{\alpha}\}, \alpha = 1, \ldots, p; \;\; i=p+1, \ldots, 
p+q$, such that for any fixed $i$, the vectors $e_i^\alpha$ 
belong to a $p$-dimensional generator $\xi$ of the Segre cone 
$SC_x (p, q)$, and for any fixed $\alpha$, the vectors 
$e_i^\alpha$ belong to a $q$-dimensional generator 
$\eta$ of  $SC_x (p, q)$.  In such a frame, the equations of 
the cone $SC_x (p, q)$ can be written as follows:
\begin{equation}\label{eq:3.1}
z_\alpha^i = t_\alpha s^i, \;\; \alpha = 1, \ldots, p, \;\; 
i = p +1, \ldots, p + q,
\end{equation}
where $z^i_\alpha$ are the coordinates of a vector 
$z = z^i_\alpha e_i^\alpha \subset T_x (M)$, and $t_\alpha$ and 
$s^i$ are parameters on which a vector $z \subset SC_x (M)$ 
depends. 

The family of frames $\{e_i^\alpha\}$ attached to the cone 
$SC_x (p, q)$ admits a transformation of the form 
\begin{equation}\label{eq:3.2}
'e_i^\alpha = A_\beta^\alpha A_i^j e_j^\beta,
\end{equation}
where $(A_\beta^\alpha)$ and $(A_j^i)$ are nonsingular square 
matrices of orders $p$ and $q$, respectively. These matrices are 
not defined uniquely since they admit a multiplication by  
reciprocal scalars. However, they can be made unique by 
restricting to unimodular matrices  $(A_\beta^\alpha)$ or 
$(A_i^j)$: $\det (A_\beta^\alpha) = 1$ or  $\det (A_i^j) = 1$. 
Thus the structural group of the almost Grassmann structure 
defined by equations (3.2), can be represented in the form 
\begin{equation}\label{eq:3.3}
G = {\bf SL} (p) \times {\bf GL}(q) \cong {\bf GL} (p) \times 
{\bf SL} (q),
\end{equation}
where ${\bf SL} (p)$ and ${\bf SL} (q)$ are special linear groups 
of dimension $p$ and $q$, respectively. Such a representation has 
been used by  Hangan [H 66, 68, 80], Goldberg [G 75a] (see also 
the book [G 88], Chapter 2), and  Mikhailov [M 78]. Unlike this 
approach, we will assume that both matrices  
$(A_{\beta}^{\alpha})$ and  $(A_i^j)$ are unimodular but 
the right-hand side of equation (3.2) admits a multiplication 
by a scalar factor. As a result, we obtain a more symmetric 
representation of the group $G$:
\begin{equation}\label{eq:3.4}
G = {\bf SL} (p) \times {\bf SL}(q)  \times {\bf H},
\end{equation}
where ${\bf H} = {\bf R}^* \otimes \mbox{{\rm Id}}$ is the group  
of homotheties of the $T_x (M)$. 

It follows that {\em an almost Grassmann structure  $AG (m, n)$ 
is a $G$-structure of first order.}

It follows from condition (3.1) that $p$-dimensional plane 
generators $\xi$ of the Segre cone $SC_x (p, q)$ are determined 
by values of the parameters $s^i$, and $t_\alpha$ are coordinates 
of points of a generator $\xi$. But a plane generator $\xi$ is 
not changed if we multiply the parameters $s^i$ by the same 
number. Thus, the family of plane generators $\xi$ depends on 
$q - 1$ parameters.

Similarly, $q$-dimensional plane generators $\eta$ of the Segre 
cone $SC_x (p, q)$ are determined by values of the parameters 
$t_\alpha$, and $s^i$ are coordinates of points of a generator 
$\eta$. 
But a plane generator $\eta$ is not changed if we multiply 
the parameters $t_\alpha$ by the same number. Thus, the family 
of plane generators $\eta$ depends on $p - 1$ parameters.

The $p$-dimensional subspaces $\xi$ form a fiber bundle on 
the manifold $M$. The base of this bundle is the manifold $M$, 
and its fiber attached to a point $x \in M$ is the set of all 
$p$-dimensional plane generators $\xi$ of the Segre cone 
$SC_x (p, q)$. The dimension of a fiber is $q - 1$, and it is 
parametrized by means of a projective space 
$P_\alpha, \; \dim P_\alpha = q - 1$. We will denote this fiber 
bundle of $p$-subspaces by $E_\alpha = (M, P_\alpha)$. 

In a similar manner, $q$-dimensional plane generators $\eta$ of 
the Segre cone $SC_x (p, q)$ form on $M$ the fiber bundle 
$E_\beta = (M, P_\beta)$ with the base $M$ and fibers of 
dimension $p - 1 = \dim P_\beta$. The fibers are $q$-dimensional 
plane generators $\eta$ of the Segre cone $SC_x (p, q)$. 
 
Consider the manifold $M_\alpha = M \times P_\alpha$ of dimension 
$pq + q - 1$. The fiber bundle $E_\alpha$ induces on $M_\alpha$ 
the distribution $\Delta_\alpha$ of plane elements $\xi_\alpha$ 
of dimension $q$. In a similar manner, on the manifold 
$M_\beta = M \times P_\beta$ the fiber bundle $E_\beta$ induces 
the distribution $\Delta_\beta$ of plane elements $\eta_\beta$ of 
dimension $p$. 

\defin{\label{def:3.6}
An almost Grassmann structure $AG (p-1, p+q-1)$ is said to be 
{\em $\alpha$-semiintegrable} if the distribution $\Delta_\alpha$ 
is  integrable on this structure.  Similarly, an almost Grassmann 
structure $AG (p-1, p+q-1)$ is said to be 
{\em $\beta$-semiintegrable} if  the distribution $\Delta_\beta$ 
is  integrable on this structure. A structure $AG (p-1, p+q-1)$ 
is called {\em integrable} if it is both $\alpha$- and 
$\beta$-semiintegrable.
} 

Integral manifolds $\widetilde{V}_\alpha$ of the distribution 
$\Delta_\alpha$ of an $\alpha$-semiintegrable almost Grassmann 
structure are of dimension $p$. They are projected on the 
original manifold $M$ in the form of a submanifold $V_\alpha$ of 
the same dimension $p$, which, at any of its points, is tangent 
to the $p$-subspace $\xi_\alpha $ of the fiber bundle $E_\alpha$. 
Through each point $x \in M$, there passes a $(q-1)$-parameter 
family of submanifolds $V_\alpha$.

Similarly, integral manifolds $\widetilde{V}_\beta$ of the 
distribution $\Delta_\beta$ of a $\beta$-semiintegrable almost 
Grassmann structure are of dimension $q$. They are projected on 
the original manifold $M$ in the form of a submanifold $V_\alpha$ 
of the same dimension $q$, which, at any of its points, is 
tangent to the $q$-subspace $\eta_\beta$ of the fiber bundle 
$E_\beta$. Through each point $x \in M$, there passes a 
$(p-1)$-parameter family of submanifolds $V_\beta$. 

If an almost Grassmann structure on $M$ is integrable, then 
through each point $x \in M$, there pass a $(q-1)$-parameter 
family of submanifolds $V_\alpha$ and a $(p-1)$-parameter family 
of submanifolds $V_\beta$ which were described above.

The Grassmann structure $G (m, n)$ is an integrable almost 
Grassmann structure $AG (m, n)$ since through 
any point $x \in \Omega (m, n)$, onto which the manifold 
$G (m, n)$ is mapped bijectively  under the Grassmann mapping, 
there pass a $(q - 1)$-parameter family of $p$-dimensional plane 
generators (which are the submanifolds $V_\alpha$) and a 
$(p - 1)$-parameter family of $q$-dimensional plane generators 
(which are the submanifolds $V_\beta$). In the projective space 
$P^n$,  to  submanifolds $V_\alpha$ there corresponds a family of 
$m$-dimensional subspaces belonging to a subspace of dimension 
$m + 1$, and  to  submanifolds $V_\beta$ there corresponds a 
family of $m$-dimensional subspaces passing through a subspace of 
dimension $m - 1$.

{\bf 3.} 
We will now write the structure equations which the forms 
$\omega_\alpha^i$ satisfy. These structure equations differ from 
equations (2.7) only by the fact that they contain an additional 
term with the product of the basis forms:
\begin{equation}\label{eq:3.5}
 d\omega_\alpha^i =  \omega_\alpha^\beta \wedge \omega_\beta^i 
+ \omega_\alpha^j  \wedge \omega_j^i
+  \kappa \wedge \omega_\alpha^i + u_{\alpha jk}^{i\beta\gamma} 
\omega_\beta^j \wedge \omega_\gamma^k,
\end{equation}
where  
$u_{\alpha jk}^{i\beta\gamma} = -u_{\alpha kj}^{i\gamma\beta}$, 
 and as earlier (see conditions (2.6)), we have 
\begin{equation}\label{eq:3.6}
\omega^\gamma_\gamma = 0,  \;\; \omega^k_k = 0. 
\end{equation}

If we prolong equations (3.5), we can see that the quantities 
$u_{\alpha jk}^{i\beta\gamma}$ as well as the quantities 
$
u_{\alpha k}^{\beta\gamma} = u_{\alpha ik}^{i\beta\gamma} 
\;\;\mbox{and} \;\;
u_{ jk}^{i\gamma} = u_{\alpha jk}^{i\alpha\gamma}
$
 form geometric objects defined in a second-order neighborhood of 
the almost Grassmann structure $AG (p - 1, p + q - 1)$. 

The following lemma can be proved (see [AG 96], Section 
{\bf 7.2}):

\begin{lemma} 
By a reduction of third-order frames of   the almost Grassmann 
structure $AG (p - 1, p + q - 1)$, the geometric objects 
$u_{\alpha k}^{\beta\gamma}$ and $u_{ jk}^{i\gamma}$ can be 
reduced to $0$:
$$
u_{\alpha k}^{\beta\gamma} = 0, \;\;
u_{ jk}^{i\gamma} = 0.
$$
\end{lemma}

The reduction indicated in Lemma 3.7 can be carried out 
by means of the fiber forms $\pi_{\alpha k}^{\beta \gamma}$ and 
$\pi_{jk}^{i \gamma}$ which appear if one finds exterior 
differentials of the forms $\omega_\alpha^\beta$ and 
$\omega_j^i$.

If we denote by $a^{i\beta\gamma}_{\alpha jk}$ the values of the 
quantities $u^{i\beta\gamma}_{\alpha jk}$ in a reduced third 
order frame, then the quantities $a^{i\beta\gamma}_{\alpha jk}$ 
satisfy the conditions
\begin{equation}\label{eq:3.7}
a^{i\alpha\gamma}_{\alpha jk} = 0, \;\;\;
a^{i\beta\gamma}_{\alpha ik} = 0,
\end{equation}
\begin{equation}\label{eq:3.8}
a^{i\beta\gamma}_{\alpha jk} = - a^{i\gamma\beta}_{\alpha kj},
\end{equation}
and 
\begin{equation}\label{eq:3.9}
\nabla_\delta a^{i\beta\gamma}_{\alpha jk} + 
a^{i\beta\gamma}_{\alpha jk} \pi = 0,
\end{equation}
where $\nabla_\delta a^{i\beta\gamma}_{\alpha jk} 
=\delta a^{i\beta\gamma}_{\alpha jk}  
- a^{i\beta\gamma}_{\sigma jk} \pi^\sigma_\alpha 
- a^{i\beta\gamma}_{\alpha sk} \pi^s_j 
- a^{i\beta\gamma}_{\alpha js} \pi^s_k 
+ a^{s\beta\gamma}_{\alpha jk} \pi_s^i 
+ a^{i\sigma\gamma}_{\alpha jk} \pi_\sigma^\beta 
+ a^{i\beta\sigma}_{\alpha jk} \pi_\sigma^\gamma, \linebreak  
\pi_m^\varepsilon = \omega_m^\varepsilon (\delta), \; 
\pi_\alpha^\beta = \omega_\alpha^\beta (\delta), \; 
\pi_i^j = \omega_i^j (\delta)$. 
This implies the following theorem.

\begin{theorem}
 The quantities $a^{i\beta\gamma}_{\alpha jk}$, defined in a 
second-order neighborhood by the reduction of third-order frames 
indicated above, form a relative tensor of weight $- 1$ 
and satisfy conditions $(3.7)$ and $(3.8)$.
\end{theorem}

\defin{\label{def:3.9}
The tensor $a = \{a_{\alpha jk}^{i\beta\gamma}\}$ is said to 
be the {\it  first structure tensor}, or the {\it torsion 
tensor}, of an almost Grassmann manifold $AG (p-1, p+q-1)$.
}

After the reduction of third-order frames has been made, the 
first structure equations (3.5) become
\begin{equation}\label{eq:3.10}
 d\omega_\alpha^i = \omega_\alpha^j \wedge 
\omega_j^i + \omega_\alpha^\beta \wedge \omega_\beta^i 
+  \kappa \wedge \omega_\alpha^i + a_{\alpha jk}^{i\beta\gamma} 
\omega_{\beta}^{j} \wedge \omega_{\gamma}^{k}.
\end{equation}

The expressions of the components of the tensor 
$a = \{a_{\alpha jk}^{i\beta\gamma}\}$ in the general (not 
reduced) third-order frame was found in [G 75] (see also [G 88], 
Section {\bf 2.2}). These expressions are:

\begin{equation}\label{eq:3.11}
\renewcommand{\arraystretch}{1.3}
\begin{array}{lll}
a_{\alpha jk}^{i\beta\gamma} & = & u_{\alpha jk}^{i\beta\gamma}-
\displaystyle \frac{2}{q^{2}-1} 
\delta_{[j}^{i}\Bigl(qu_{|\alpha|k]}^{[\beta\gamma]}+
 u_{|\alpha|k]}^{[\gamma\beta]}\Bigr)  
  -\displaystyle \frac{2}{p^{2}-1}\delta_{\alpha}^{[\beta}
\Bigl(pu_{[jk]}^{|i|\gamma]}+
 u_{[kj]}^{|i|\gamma]}\Bigr)  \\
      &   & + \displaystyle \frac{2}{(p^{2}-1)(q^{2}-1)}
 \Bigl[(pq-1)\Bigl(\delta_{[j}^{i}\delta_{|\alpha|}^{[\beta}
 u_{k]}^{\gamma]} 
+ \delta_{[k}^i\delta_{|\alpha|}^{[\beta}
 \widetilde{u}_{j]}^{\gamma]}\Bigr)  \\
      &  & 
+ (q-p)\Bigl(\delta_{[j}^i\delta_{|\alpha|}^{[\beta} 
 \widetilde{u}_{k]}^{\gamma]}
+ \delta_{[k}^i\delta_{|\alpha|}^{[\beta} 
u_{j]}^{\gamma]}\Bigr) \Bigr],
\end{array}
\renewcommand{\arraystretch}{1}
\end{equation}
where the alternation is carried out with respect 
to the pairs of indices ${\beta \choose j}$, ${\gamma \choose k}$ 
or ${\beta \choose k}$, ${\gamma \choose j}$, and 
$
u_k^\gamma = u_{lk}^{l\gamma} = u_{\sigma lk}^{l\sigma\gamma} 
= u_{\sigma k}^{\sigma\gamma}, \;\;
\widetilde{u}_k^\gamma = u_{kl}^{l\gamma} 
= u_{\sigma kl}^{l\sigma\gamma} = - u_{\sigma lk}^{l\gamma\sigma} 
= - u_{\sigma k}^{\gamma\sigma}.
$

If we prolong equations (3.10) and make a reduction of fourth 
order frames of the almost Grassmann structure 
$AG (p - 1, p + q - 1)$, we will find the following remaining 
structure equations of $AG (p - 1, p + q - 1)$ which the forms 
$\omega_\alpha^\beta, \omega_j^i$ and $\kappa$ satisfy: 
\begin{equation}\label{eq:3.12}
\renewcommand{\arraystretch}{1.3}
 \begin{array}{ll}
d \omega_\alpha^\beta - \omega_\alpha^\gamma \wedge \omega_\gamma^\beta 
= \displaystyle \frac{q}{p+q}\Bigl(\delta_\alpha^\beta 
\omega_\gamma^k \wedge \omega^\gamma_k - p \omega_\alpha^k 
\wedge \omega_k^\beta\Bigr) 
+  b_{\alpha kl}^{\beta\gamma\delta}  \omega_\gamma^k \wedge 
\omega_\delta^l, \\
 d \omega_j^i - \omega_j^k  \wedge \omega_k^i
= \displaystyle \frac{p}{p+q} \Bigl(\delta_j^i \omega^\gamma_k 
\wedge \omega_\gamma^k 
- q \omega_j^\gamma  \wedge \omega_\gamma^i\Bigr) 
+  b_{jkl}^{i\gamma\delta}  \omega_\gamma^k \wedge 
\omega_\delta^l, \\
d \kappa = \omega_i^\alpha \wedge \omega_\alpha^i,
\end{array}
\renewcommand{\arraystretch}{1}
\end{equation}
where the quantities $b_{\alpha lm}^{\beta \delta \varepsilon}$ 
and $b_{jlm}^{i \delta \varepsilon}$ are defined in a third-order 
neighborhood and  satisfy the conditions 
\begin{equation}\label{eq:3.13}
\renewcommand{\arraystretch}{1.3}
\begin{array}{ll}
b_{\alpha lm}^{\beta \delta \varepsilon} = 
- b_{\alpha ml}^{\beta  \varepsilon\delta}, & 
b_{jlm}^{i \delta \varepsilon} = 
- b_{jml}^{i  \varepsilon\delta}, \\
b_{\sigma lm}^{\sigma \delta \varepsilon} = 0, &
b_{klm}^{k \delta \varepsilon} = 0, \\  b^{\gamma\alpha\delta}_{\alpha kl} 
-  b^{i\gamma\delta}_{kil}  
+ b^{\delta\alpha\gamma}_{\alpha lk} 
- b^{i\delta\gamma}_{lik} = 0. &
\end{array}
\renewcommand{\arraystretch}{1}
\end{equation}
Note that the last  conditions in (3.12) and (3.13) are 
the result of reduction of fourth-order frames mentioned above.  

If we prolong equations (3.12) and set $\omega^i = 0$ in the 
resulting equations, we find that the quantities  
$b_{\alpha lm}^{\beta \delta \varepsilon}$ and 
$b_{jlm}^{i \delta \varepsilon}$ satisfy the following equations: 
\begin{equation}\label{eq:3.14}
\renewcommand{\arraystretch}{1.3}
\begin{array}{ll}
\nabla_\delta b^{\beta\gamma\delta}_{\alpha kl} 
+ 2 b^{\beta\gamma\delta}_{\alpha kl} \pi 
- \displaystyle \frac{pq}{2(p+q)} 
\Bigl(2\delta_\varepsilon^\beta a^{m\gamma\delta}_{\alpha kl} 
- \delta_\alpha^\gamma a^{m\beta\delta}_{\varepsilon kl} 
+ \delta_\alpha^\delta a^{m\beta\gamma}_{\varepsilon lk}\Bigr)  
\pi_m^\varepsilon = 0, \\
\nabla_\delta b^{i\gamma\delta}_{jkl} 
+ 2 b^{i\gamma\delta}_{jkl} \pi 
+ \displaystyle \frac{pq}{2(p+q)} 
\Bigl(2\delta_j^m a^{i\gamma\delta}_{\varepsilon kl}
- \delta_k^i a^{m\gamma\delta}_{\varepsilon jl} 
+ \delta_l^i a^{i\delta\gamma}_{\varepsilon jk}\Bigl)
 \pi_m^\varepsilon = 0,
\end{array}
\renewcommand{\arraystretch}{1}
\end{equation}
where the operator $\nabla_\delta$ is defined in the same way as 
in formula (3.9), and 
$\delta$ is the symbol of differentiation with respect to the 
fiber parameters. It follows from equations (3.14) that the 
quantities $\{b_{\alpha lm}^{\beta \gamma \delta}\}$ and 
$\{b_{jlm}^{i \alpha \beta}\}$ do not form tensors or even 
homogeneous geometric objects but the quantities 
$\{b_{\alpha lm}^{\beta \gamma \delta}, \; a_{\alpha jk}^{i\beta 
\gamma} \}$ as well as the quantities  
$\{b_{jlm}^{i \alpha \beta}, \; a_{\alpha jk}^{i\beta \gamma}\}$  
form linear homogeneous objects. They represent two subobjects of 
the {\em second structure object}  (or the {\em torsion-curvature 
object}) $\{a_{\alpha jk}^{i\beta \gamma}, \; 
b_{\alpha lm}^{\beta \gamma \delta},\; 
b_{jlm}^{i \alpha \beta}\}$  
of the almost Grassmann structure \\$AG (p-1, p+q-1)$.

Moreover, the prolongation of equations (3.12) leads to the 
following structure equations:
\begin{equation}\label{eq:3.15}
d \omega_i^\alpha - \omega_i^\beta \wedge \omega_\beta^\alpha 
- \omega_i^j \wedge \omega_j^\alpha  
+ \kappa \wedge \omega_i^\alpha
=    c_{ijk}^{\alpha\beta\gamma} 
 \omega^k_\gamma \wedge \omega_\beta^j 
- a_{\gamma ij}^{k\alpha\beta} \omega_k^\gamma \wedge 
\omega_\beta^j,
\end{equation}
where $c_{ijk}^{\alpha\beta\gamma} 
= - c_{ikj}^{\alpha\gamma\beta}$. 

In addition, taking exterior derivatives of the last equation of 
(3.12) and applying (3.10) and (3.15), we find the following 
condition for the quantities $c_{ijk}^{\alpha\beta\gamma}$:
\begin{equation}\label{eq:3.16}
c_{[ijk]}^{[\alpha\beta\gamma]} = 0.
\end{equation}

Finally, if we prolong equations (3.15), we find that 
\begin{equation}\label{eq:3.17}
\renewcommand{\arraystretch}{1.3}
\begin{array}{ll}
\bigl[\nabla c_{ijk}^{\alpha \beta \gamma}  
&+ 3 c_{ijk}^{\alpha \beta  \gamma}  \kappa 
- b_{\sigma kj}^{\alpha \gamma \beta} \omega_i^\sigma 
+ b_{ikj}^{l \gamma \beta} \omega_l^\alpha 
+ (a^{l \alpha\beta}_{\delta ij} 
a^{m\delta\gamma}_{\varepsilon lk} 
+ a^{m \alpha\delta}_{\varepsilon il} 
a^{l\gamma\beta}_{\delta kj}) \omega^\varepsilon_m \\
&+ (2 c_{iml}^{\alpha\varepsilon\delta} 
a_{\beta kj}^{l\gamma\beta}  \omega_\varepsilon^m 
- c_{ljk}^{\delta\beta\gamma} 
a_{\delta im}^{l \alpha\varepsilon}) 
\omega_\varepsilon^m\bigr] \wedge \omega^k_\gamma \wedge \omega^j_\beta = 0,
\end{array}
\renewcommand{\arraystretch}{1}
\end{equation}
where $\nabla c_{ijk}^{\alpha \beta \gamma} 
= d c_{ijk}^{\alpha \beta \gamma}  
-  c_{ljk}^{\alpha \beta \gamma}  \omega_i^l 
-  c_{ilk}^{\alpha \beta \gamma}  \omega_j^l 
 -  c_{ijl}^{\alpha \beta \gamma}  \omega_k^l 
+  c_{ijk}^{\delta \beta \gamma}  \omega_\delta^\alpha 
+  c_{ijk}^{\alpha \delta  \gamma}  \omega_\delta^\beta 
+  c_{ijk}^{\alpha  \beta \delta}  \omega_\delta^\gamma$. 
 For $\omega_\alpha^i = 0$, it follows from equation (3.17) that  
\begin{equation}\label{eq:3.18}
\nabla_\delta c_{ilk}^{\alpha \delta \gamma}  
+ 3 c_{ilk}^{\alpha \beta \gamma}  \pi 
- b_{\sigma kl}^{\alpha \gamma \delta} \pi_i^\sigma 
+ b_{ikl}^{j  \gamma \delta} \pi_j^\alpha 
+ (a_{\beta il}^{j\alpha\delta} a_{\varepsilon jk}^{m\beta\gamma} 
+ a_{\varepsilon ij}^{m\alpha\beta} a_{\beta kl}^{j\gamma\delta}) 
\pi_m^\varepsilon = 0.
\end{equation}
Equations  (3.9), (3.13) and (3.18) prove that the quantities 
$\{a_{\alpha jk}^{i\beta \gamma},\; 
b_{\alpha km}^{\beta \gamma \delta}, \; 
b_{ijk}^{l \beta \gamma}, \; c_{ijk}^{\alpha \beta \gamma} \}$ 
also form a linear homogeneous object which is called the {\em 
third structure object} of the almost Grassmann structure 
$AG (m, n)$. Note that the torsion tensor 
$\{a_{\alpha jk}^{i\beta \gamma}\}$ is defined in a 
second-order 
differential neighborhood of a point $x \in M$, the second 
structure object  $ \{a_{\alpha jk}^{i\beta \gamma}, \; 
b_{\alpha km}^{\beta \gamma \delta},\; 
b_{jlm}^{i \alpha \beta}\}$ is defined in a third-order 
differential 
neighborhood of a point $x \in M$, and the third structure object 
$\{a_{\alpha jk}^{i\beta \gamma},\; b_{\alpha km}^{\beta \gamma 
\delta}, \; b_{ijk}^{l  \beta \gamma}, \; 
c_{ijk}^{\alpha \beta \gamma}\}$ is defined in a fourth 
differential neighborhood of a point $x \in M$. 

The structure equations (3.10), (3.12) and (3.15) can be 
written in the index-free notation as follows: 
\begin{equation}\label{eq:3.19}
\renewcommand{\arraystretch}{1.3}
\left\{
\begin{array}{ll}
 d \omega =  \kappa \wedge \omega -  \omega  \wedge \theta_\alpha 
-  \theta_\beta \wedge \omega + \Omega, \\
d \theta_\alpha  + \theta_\alpha \wedge \theta_\alpha  
= \displaystyle \frac{q}{p+q} \Bigl[- I_\alpha 
\mbox{{\rm tr}}\; (\varphi \wedge \omega) 
+ p \varphi \wedge \omega\Bigr] 
+ \Theta_\alpha,  \\
d \theta_\beta  + \theta_\beta \wedge \theta_\beta  
= \displaystyle \frac{p}{p+q} \Bigl[- I_\beta 
\mbox{{\rm tr}}\; (\varphi \wedge \omega) 
+ q \omega \wedge \varphi\Bigr] 
+ \Theta_\beta,  \\
d \kappa  =  \mbox{{\rm tr}} \; (\varphi \wedge \omega), \\
 d \varphi  + \theta_\alpha \wedge \varphi  + \varphi \wedge 
\theta_\beta  + \kappa \wedge  \varphi = - (a \varphi) 
\wedge \omega + \Phi,
\end{array}
\right.
\renewcommand{\arraystretch}{1}
\end{equation}
where $\omega = (\omega_\alpha^i)$ is a 
matrix 1-form defined in the first order frame bundle;  
$\kappa$ is a  scalar 1-form,  
$\theta_\alpha = (\omega_\beta^\alpha)$ and 
$\theta_\beta = (\omega_j^i)$ are the matrix 1-forms 
defined in the second-order frame bundle for which 
$$
\mbox{{\rm tr}} \; \theta_\alpha = 0, \;\; 
\mbox{{ \rm tr}} \; \theta_\beta = 0;
$$
$\varphi = (\omega_i^\alpha)$ is a matrix 1-form 
defined in a third-order fiber bundle; 
 $I_\alpha = (\delta_\alpha^\beta)$ and 
$I_\beta = (\delta_i^j)$ 
are the unit tensors of orders $p$ and $q$, 
respectively; the 2-form 
$\Omega = (\Omega_\alpha^i)$ is the torsion form; and the 2-forms 
  $\Theta_\alpha = (\Theta^\alpha_\beta), \;  \Theta_\beta 
= (\Theta^i_j)$, and $\Phi = (\Phi^\alpha_i)$ are the 
 {\em curvature $2$-forms} of the  $AG (p-1, p+q-1)$-structure. 
The  components of 2-forms $\Omega, \; \Theta_\alpha, \;  
\Theta_\beta$, and $\Phi$ are
\begin{equation}\label{eq:3.20}
\renewcommand{\arraystretch}{1.3}
\begin{array}{ll} 
 \Omega_\alpha^i =  a^{i\beta\gamma}_{\alpha jk} \omega_\gamma^k 
 \wedge \omega_\delta^l, &  
\Theta^\alpha_\beta = b_{\beta kl}^{\alpha \gamma \delta} \omega_\gamma^k \wedge \omega_\delta^l, \\ 
\Theta^i_j = b_{jkl}^{i \gamma \delta} \omega_\gamma^k \wedge 
\omega_\delta^l, & 
\Phi_i^\alpha = c_{ikl}^{\alpha \gamma \delta} 
\omega_\delta^l \wedge \omega_\gamma^k.
\end{array}
\renewcommand{\arraystretch}{1}
\end{equation}

{\bf 4.} The restrictions of equations (3.19)  to  a fiber  
frame bundle, i.e. for  $\omega = 0$, have the form
\begin{equation}\label{eq:3.21}
\renewcommand{\arraystretch}{1.3}
\begin{array}{ll}
d \kappa = 0,  \\
d \theta_\alpha  = - \theta_\alpha \wedge \theta_\alpha, \;\; 
d \theta_\beta = - \theta_\beta \wedge \theta_\beta, \\
d \varphi = -\kappa \wedge \varphi - \theta_\alpha \wedge 
\varphi - \varphi \wedge \theta_\beta. 
\end{array}
\renewcommand{\arraystretch}{1}
\end{equation}

From equations (3.21) it follows that 

\begin{description}

\item[(1)] The form $\kappa$ is an invariant form of the group 
${\bf H}$ of homotheties acting in the space $T_x (M)$. 

\item[(2)]  The forms $\theta_\alpha$ and $\theta_\beta$ are  
invariant forms of the special linear groups  ${\bf SL} (p)$ and 
${\bf SL} (q)$, respectively.

\item[(3)] The 1-forms $\kappa, \theta_\alpha$, and 
$\theta_\beta$ are invariant forms of the structural group $G$ of 
the almost Grassmann structure $AG (m, n)$ whose transformations 
leave cone $SC_x (p, q)$ invariant. As noted earlier, the group 
$G$ is isomorphic to the direct product ${\bf SL} (p) \times 
{\bf SL} (q) \times {\bf H}$: 
\begin{equation}\label{eq:3.22}
G \cong {\bf SL} (p) \times {\bf SL} (q) \times {\bf H}.
\end{equation}

\item[(4)] 
The 1-forms  $\kappa, \theta_\alpha, \theta_\beta$. and $\varphi$ 
are invariant forms of the group $G'$ which is a differential 
prolongation of $G$. The group $G'$ is isomorphic to the group 
$G \htimes {\bf T} (pq)$ whose subgroup $ {\bf T} (pq)$ is 
defined by the invariant forms  $\omega_i^\alpha$. Thus, 
\begin{equation}\label{eq:3.23}
G' \cong  ({\bf SL} (p) \times {\bf SL} (q) \times 
{\bf H}) \htimes {\bf T} (pq).
\end{equation}
\end{description}

 Since the group $G'$ does not admit further prolongations,  
{\em an  almost Grassmann structure $AG (m, n)$ is a 
$G$-structure of finite type two}. 

In order to describe the group $G'$ geometrically, we compactify 
the tangent subspace $T_x (M)$ by enlarging it by the point 
at infinity, $y = \infty$,  and the Segre cone $SC_y (p, q)$ 
with its vertex at this point. Then the manifold 
$T_x (M) \cap SC_y (p, q)$ is equivalent to the algebraic 
variety $\Omega (p, q)$. Since the point $x$, at which the 
variety $\Omega (p, q)$ is tangent to the manifold $M$, is fixed, 
the geometry defined by the group $G'$ on $\Omega (p, q)$ is 
equivalent to that of the space obtained as a projection of the 
variety $\Omega (p, q)$ from the point $x$  onto a flat space of 
dimension $pq$. The group $G'$ is the group of motions of this 
space,  subgroup $G$ is the isotropy group of this space, and 
the subgroup ${\bf T} (pq)$ is the subgroup of parallel 
translations.

The group $G'$ can be also represented as the group of motions of 
a projective space $P^n$ leaving invariant a fixed subspace 
$P^m = x$. The subgroup ${\bf SL} (p)$ of $G'$ is locally 
isomorphic to the group of projective transformations of the 
subspace $P^m$; the subgroup ${\bf SL} (q)$ of $G'$ is 
locally isomorphic to the group of transformations in 
the bundle of $(m+1)$-dimensional subspaces of $P^n$ passing 
through $x$; its the subgroup  ${\bf H}$ of $G'$ 
is the group of homotheties 
with its center at $x$; and the subgroup ${\bf T} (pq)$ of $G'$ 
is the group of translations of the  subspaces $y = P^{n-m-1}$ 
 which are complementary to $x = P^m$ in $P^n$.

{\bf 5.} As it was proved in [G 75] (see also [G 88] and [M 78]), 
the torsion tensor $a = \{a_{\alpha jk}^{i\beta\gamma}\}$ of an 
almost Grassmann structure $AG (p-1, p+q-1)$ satisfying 
conditions (3.7)--(3.8) decomposes into two subtensors:
\begin{equation}\label{eq:3.24}
a = a_\alpha \dot{+} a_\beta, 
\end{equation}
where
$$
a_\alpha = \{a_{\alpha (jk)}^{i\beta\gamma}\}  \;\;\;
\mbox{and}\;\;\; 
a_\beta = \{a_{\alpha jk}^{i(\beta\gamma)}\}.   
 $$
Note that $a_{\alpha (jk)}^{i\beta\gamma} 
= a_{\alpha jk}^{i[\beta\gamma]}$ and 
$
a_{\alpha jk}^{i(\beta\gamma)} = a_{\alpha [jk]}^{i\beta\gamma}.
$
It is easy to see that  the components of subtensors $a_\alpha$ 
and $a_\beta$ satisfy the conditions similar to conditions 
 (3.8).

In addition, it is easy to prove that:

\begin{description}
\item[(i)]  {\it If $p=2$, then $a_\alpha = 0$}. 
\item[(ii)] {\it If $q=2$, then $a_\beta = 0$}. 
\end{description}

 There are certain dependencies between three structure 
objects indicated above. Let us set 
$b^1 = \{b^{i\gamma\delta}_{jkl}\}, \;\; b^2 = 
 \{b^{\alpha\gamma\delta}_{\beta kl}\}$,  
and $\{c^{\alpha\gamma\delta}_{ikl}\}$. 
 It can be proved  that 

\begin{description}
\item[(1)] {\it If $q > 2$, then the components of  $b^2$  
are expressed in terms of components of the tensor $a$ and their 
Pfaffian derivatives, and the components of  $c$ are expressed in 
terms of components of the object $(a, b^2)$ and their Pfaffian  
derivatives.}

\item[(2)] {\it If $p > 2$, then the components of $b^1$ are 
expressed in terms of components of 
the tensor $a$ and their Pfaffian 
derivatives, and the components of  $c$ are expressed in terms of 
components of the  $(a, b^1)$  and their Pfaffian derivatives}. 

\item[(iii)] {\it If $p > 2$ and $q > 2$, then the components of  
$b$ and $c$ are expressed in terms of components of the tensor 
$a$ and their Pfaffian  derivatives.} 

\end{description}

However, the tensor $a$ itself is not arbitrary since if we 
substitute for the components $b$ and $c$ their expressions in 
terms of the tensor $a$ and its Pfaffian derivatives into 
conditions of integrability of the structure equations (3.12) and 
(3.15), we obtain certain algebraic conditions for the components 
of the tensor $a$ and its covariant derivatives. The latter 
conditions are analogues of the Bianchi equations in the theory 
of spaces with affine connection.

The structure object $S = \{a, b, c\}$ of the almost Grassmann 
structure $AG (p-1, p+q-1)$ is {\em complete} in the sense that 
if we
 prolong the structure equations (3.19) of $AG (p-1, p+q-1)$, 
then all newly arising objects are expressed in terms of the 
components of the object $S$ and their Pfaffian derivatives. 
This follows from the fact that {\em the almost Grassmann 
structure 
$AG (p-1, p+q-1)$ is a $G$-structure of finite type two}. 

\defin{\label{def:2.8} 
 An almost  Grassmann structure $AG (p-1, p+q-1)$ is said to be 
{\em locally Grassmann} (or {\em locally flat}) if it is locally 
equivalent to a Grassmann structure. }

This means that a locally flat almost  Grassmann structure 
$AG (p-1, p+q-1)$ admits a mapping onto an open domain of the 
algebraic variety $\Omega (m, n)$ of a projective space 
$P^N$, where $N = \displaystyle {n+1 \choose m+1} - 1, \; 
m = p - 1, n = p + q - 1$, under which the Segre cones of 
the  structure $AG (p-1, p+q-1)$ correspond to the asymptotic 
cones of variety $\Omega (m, n)$.

 From the equivalence theorem of \'{E}. Cartan (see [C 08] or  
[Ga 89]), it follows that in order for an almost  Grassmann 
structure $AG (p-1, p+q-1)$ to be locally Grassmann, it is 
necessary and sufficient that its structure equations have the 
form (2.17). Comparing these equations with equations (3.19), we 
see that  {\em an almost  Grassmann structure $AG (p-1, p+q-1)$ 
is locally Grassmann if and only if its complete structure object 
$S = (a, b, c)$ vanishes.}

However, as was noted above, if $p > 2$ and $q > 2$, the 
components of $b$ are expressed in terms of the components of the 
tensor $a$ and their Pfaffian derivatives, and the components of 
$c$ are expressed in terms of the components of the subobject 
$(a, b)$ and their Pfaffian derivatives. Moreover, it can be 
proved that the vanishing of the tensor $a$ on a manifold $M$ 
carrying an almost  Grassmann structure implies the vanishing of 
the components of $b$ and $c$. 

This implies the  following result.

\begin{theorem} For $p > 2$ and $q > 2$,  an almost  Grassmann 
structure $AG (p-1, p+q-1)$ is locally Grassmann if and only if 
its first structure tensor $a$ vanishes. 
\end{theorem}

\section{Semiintegrability of Almost Grassmann Structures}

\setcounter{equation}{0}
 
{\bf 1.} 
Now we will prove the following necessary and sufficient 
conditions for the almost Grassmann structure $AG (p-1, p+q-1)$ 
to be $\alpha$- or  $\beta$-semiintegrable.

\begin{theorem}

\begin{description}
\item[(1)] 
If $p > 2$ and $q \geq 2$, then for an almost Grassmann structure 
$AG (p-1, p+q-1)$ to be $\alpha$-semiintegrable, it is necessary 
and sufficient that the following condition holds: 
$a_\alpha = b_\alpha^1 = b_\alpha^2 = 0$.

\item[(2)] 
If $p \geq 2$ and $q > 2$, then for an almost Grassmann structure 
$AG (p-1, p+q-1)$ to be $\beta$-semiintegrable, it is necessary 
and sufficient that the following conditions hold: 
$a_\beta = b_\beta^1 = b_\beta^2 = 0$. 

\end{description}
\end{theorem}

{\sf Proof.} We will prove  part (1) of theorem. The proof 
of part (2) is similar.

Suppose that $\theta_\alpha, \; \alpha = 1, \ldots, p$, are 
basis forms of the subvarieties $V_\alpha, \; \dim V_\alpha 
= p$, indicated in Definition 
{\bf 3.6}. Then 
\begin{equation}\label{eq:4.1}
\omega_\alpha^i = s^i \theta_\alpha, \;\;
 \alpha = 1, \ldots, p; \;\; i = p+1, \ldots , p+q.
\end{equation}

For the structure $AG (p-1, p+q-1)$ to be 
$\alpha$-semiintegrable, it is necessary and sufficient that 
system (4.1) be completely integrable. Taking the exterior 
derivatives of equations (4.1) by means of structure equations 
(3.10), we find that 
\begin{equation}\label{eq:4.2}
(d s^i + s^j \omega_j^i - s^i \omega) \wedge \theta_\alpha 
+ s^i (d \theta_\alpha - \omega_\alpha^\beta \wedge \theta_\beta) 
= a_{\alpha jk}^{i \beta\gamma} s^j s^k \theta_\beta \wedge 
\theta_\gamma.
\end{equation}
It follows from these equations that 
\begin{equation}\label{eq:4.3}
d \theta_\alpha -  \omega_\alpha^\beta \wedge \theta_\beta 
= \varphi_\alpha^\beta  \wedge \theta_\beta,
\end{equation}
where $\varphi_\alpha^\beta$ is an 1-form that is not 
expressed in terms of  the 
basis forms $\theta_\alpha$. 

For brevity, we set 
\begin{equation}\label{eq:4.4}
 \varphi^i = d s^i + s^j \omega_j^i - s^i \omega.
\end{equation}
Then the exterior quadratic equation (4.2) takes the form:
\begin{equation}\label{eq:4.5}
(\delta_\alpha^\beta \varphi^i + s^i \varphi_\alpha^\beta)  
\wedge \theta_\beta 
= a_{\alpha jk}^{i \beta\gamma} s^j s^k \theta_\beta \wedge 
\theta_\gamma.
\end{equation}
 From (4.5) it follows that for $\theta_\alpha = 0$, the 1-form 
$\delta_\alpha^\beta \varphi^i + s^i \varphi_\alpha^\beta$ 
vanishes:
\begin{equation}\label{eq:4.6}
\delta_\alpha^\beta \varphi^i (\delta) 
+ s^i \varphi_\alpha^\beta (\delta)= 0.
\end{equation}
Contracting equation (4.6) with respect to the indices 
$\alpha$ and $\beta$, we find that 
\begin{equation}\label{eq:4.7}
\varphi^i = - s^i \varphi (\delta), \;\; 
 \varphi_\alpha^\beta = \delta_\alpha^\beta \varphi (\delta),
\end{equation}
where we set $\varphi (\delta) =  \frac{1}{p} 
\varphi_\gamma^\gamma$.

It follows from (4.7) that on the subvariety $V_\alpha$,  
the 1-forms  
$\varphi^i$ and $\varphi_\alpha^\beta$ can be written as follows:
\begin{equation}\label{eq:4.8}
\varphi^i = - s^i \varphi + s^{i\beta} \theta_\beta, \;\; 
 \varphi_\alpha^\beta = \delta_\alpha^\beta \varphi  
+ \widehat{s}_\alpha^{\beta\gamma} \theta_\gamma.
\end{equation}
Substituting these expressions into equations (4.3) and (4.4), 
we find that 
\begin{equation}\label{eq:4.9}
d \theta_\alpha -  \omega_\alpha^\beta \wedge \theta_\beta 
=  s_\alpha^{\beta\gamma} \theta_\gamma \wedge \theta_\beta
\end{equation}
where 
$s_\alpha^{\beta\gamma} = \widehat{s}_\alpha^{[\beta\gamma]}$, 
and 
\begin{equation}\label{eq:4.10}
  d s^i + s^j \omega_j^i - s^i \kappa 
=  - s^i \varphi  + s^{i\beta} \theta_\beta.
\end{equation}
Substituting (4.9) and (4.10) into equation (4.2), we obtain
\begin{equation}\label{eq:4.11}
- s^i s_\alpha^{\beta\gamma} 
- \delta_\alpha^{[\beta} s^{|i|\gamma]} 
= a_{\alpha jk}^{i [\beta\gamma]} s^j s^k.
\end{equation}
Contracting equation (4.11) with respect to the indices $\alpha$ 
and $\beta$, we find that
$$
- 2 s^i s_\alpha^{\alpha\gamma} - p s^{i\gamma} + 
s^{i\gamma} = 0.
$$
It follows that 
\begin{equation}\label{eq:4.12}
 s^{i\gamma} = s^i s^\gamma,
\end{equation}
where we set 
$s^\gamma = - \frac{2}{p-1} s_\alpha^{\alpha\gamma}$.
Substituting (4.12) into (4.11), we find that 
\begin{equation}\label{eq:4.13}
s^i (\delta_\alpha^\gamma s^\beta 
- \delta_\alpha^\beta s^\gamma - 2 s_\alpha^{\beta\gamma}) 
= 2 a_{\alpha jk}^{i [\beta\gamma]} s^j s^k. 
\end{equation}
It follows that 
\begin{equation}\label{eq:4.14}
\delta_\alpha^\gamma s^\beta 
- \delta_\alpha^\beta s^\gamma - 2 s_\alpha^{\beta\gamma} 
=  s_{\alpha j}^{\beta\gamma} s^j, 
\end{equation}
where $s_{\alpha j}^{\beta\gamma} 
= - s_{\alpha j}^{\gamma\beta}$. 
Substituting (4.14) into (4.13), we arrive at the equation 
\begin{equation}\label{eq:4.15}
 s_{\alpha (j}^{\beta\gamma} \delta^i_{k)} = 
a_{\alpha (jk)}^{i\beta\gamma}, 
\end{equation}
where the alternation sign in the right-hand side was dropped 
since $a_{\alpha (jk)}^{i\beta\gamma} 
= a_{\alpha jk}^{i[\beta\gamma]}$.

Contracting (4.15) with respect to the indices $i$ and $j$ 
and taking into account equations (3.7) and (3.8), we obtain 
\begin{equation}\label{eq:4.16}
s_{\alpha k}^{\beta\gamma} = 0. 
\end{equation}
By (4.15), it follows that 
\begin{equation}\label{eq:4.17}
a_{\alpha (jk)}^{i\beta\gamma} = 0. 
\end{equation}

Thus, we proved that {\em if an almost Grassmann structure 
$AG (p - 1, p + q - 1)$ is $\alpha$-semiintegrable, then its 
torsion tensor satisfies the condition $(4.17)$, i.e. 
$a_\alpha = 0$.} 

Since as was noted earlier, for $p = 2$, the subtensor 
$a_\alpha = 0$,  condition (4.17) is identically satisfied. 
Hence while proving sufficiency of this condition for 
$\alpha$-semiintegrability, we must assume that $p > 2$.

Let us return to equations (4.9) and (4.10). Substitute 
into equation (4.10) the values $s^{i\beta}$ taken from 
(4.12) and set 
\begin{equation}\label{eq:4.18}
\widetilde{\varphi} = \varphi - s^\beta \theta_\beta.
\end{equation}
In addition, by (4.16), relations (4.14) imply that
$$
s_\alpha^{\beta\gamma} = \delta_\alpha^{[\gamma} s^{\beta]}.
$$
Then equations (4.9) and (4.10) take the form 
\begin{equation}\label{eq:4.19}
d \theta_\alpha -  (\omega_\alpha^\beta 
+ \delta_\alpha^\beta \widetilde{\varphi}) \wedge \theta_\beta
= 0 
\end{equation}
and 
\begin{equation}\label{eq:4.20}
d s^i + s^j \omega_j^i - s^i (\kappa - \widetilde{\varphi}) = 0.
\end{equation}
Taking the exterior derivatives of (4.20), we obtain the 
following exterior quadratic equation: 
\begin{equation}\label{eq:4.21}
s^i \Phi + b_{jkl}^{i\gamma\delta} s^j s^k s^l \theta_\gamma 
\wedge \theta_\delta = 0,
\end{equation}
where 
$$
\Phi = d \widetilde{\varphi} - \displaystyle \frac{(p+1)q}{p+q} 
s^k \omega_k^\gamma \wedge \theta_\gamma.
$$
Next, taking the exterior derivatives of (4.19), we find that 
\begin{equation}\label{eq:4.22}
 \Phi \wedge \theta_\alpha + b_{\alpha kl}^{\beta\gamma\delta}
 s^k s^l \theta_\beta \wedge \theta_\gamma \wedge \theta_\delta 
= 0.
\end{equation}
Equation (4.21) shows that the 2-form $\Phi$ 
can be written as 
\begin{equation}\label{eq:4.23}
 \Phi = A_{kl}^{\gamma\delta} s^k s^l  \theta_\gamma \wedge \theta_\delta,
\end{equation}
where the coefficients $A_{kl}^{\gamma\delta}$ are symmetric 
with respect to the lower indices and skew-symmetric 
with respect to the upper indices. Substituting this value of the 
form $\Phi$ into equations (4.21) and (4.22), we arrive 
at the conditions:
\begin{equation}\label{eq:4.24}
b_{(jkl)}^{i[\gamma\delta]} 
+ \delta^i_{(j} A_{kl)}^{\gamma\delta} = 0
\end{equation}
and 
\begin{equation}\label{eq:4.25}
b_{\alpha (kl)}^{[\beta\gamma\delta]}
 + \delta_\alpha^{[\beta} A_{kl}^{\gamma\delta]} = 0.
\end{equation}
Contracting equation (4.24) with respect to the indices 
$i$ and $j$ and  equation (4.25) with respect to the indices 
$\alpha$ and $\beta$, we obtain 
\begin{equation}\label{eq:4.26}
2 (q + 2) A_{kl}^{\gamma\delta} 
+ b_{kli}^{i\gamma\delta} + b_{kil}^{i\gamma\delta} 
+  b_{lik}^{i\gamma\delta} + b_{lki}^{i\gamma\delta}  = 0
\end{equation}
and 
\begin{equation}\label{eq:4.27}
2 (p - 2) A_{kl}^{\gamma\delta} 
+ b_{\alpha kl}^{\gamma\delta\alpha} 
+ b_{\alpha lk}^{\gamma\delta\alpha} 
+ b_{\alpha kl}^{\delta\alpha\gamma} 
+ b_{\alpha lk}^{\delta\alpha\gamma} = 0.
\end{equation}
Note that for $p = 2$,, equation (4.25) becomes an identity, and 
we will not obtain equations (4.27). 

If we add equations (4.26) and (4.27) and apply the last 
condition of (3.13), we find that 
\begin{equation}\label{eq:4.28}
A_{kl}^{\gamma\delta} = 0.
\end{equation}
As a result, equations (4.24) and (4.25) take the form
\begin{equation}\label{eq:4.29}
b_{(jkl)}^{i[\gamma\delta]} = 0, \;\;\; 
b_{\alpha (kl)}^{[\beta\gamma\delta]} = 0.
\end{equation}

By the first conditions of (3.13), conditions (4.29) are 
equivalent to the conditions
\begin{equation}\label{eq:4.30}
b_{(jkl)}^{i\gamma\delta} = 0, \;\;\; 
b_{\alpha kl}^{[\beta\gamma\delta]} = 0.
\end{equation}

It follows from equations (4.28) and (4.23) that 
\begin{equation}\label{eq:4.31}
d \widetilde{\varphi} = \displaystyle \frac{(p+1)q}{p+q} 
s^k \omega_k^\gamma \wedge \theta_\gamma.
\end{equation}
Finally, taking the exterior derivatives of equations (4.31) and 
applying (4.19), (4.20) and (3.15), we obtain the condition
\begin{equation}\label{eq:4.32}
 c_{(ijk)}^{[\alpha\beta\gamma]} = 0.
\end{equation}
This equation will not be trivial only if $p > 2$. It is easy to 
check that the last condition follows from integrability 
condition (3.16).

Thus, the system of Pfaffian equations (4.1), defining integral 
submanifolds of an $\alpha$-semiintegrable almost Grassmann 
structure, together with Pfaffian equations (4.10) and (4.31) 
following from (4.1) is completely integrable if and only if 
conditions (4.17) and (4.30) are satisfied. This 
concludes the proof of part (I) of the theorem. 

We introduce the following notations:
$$
\begin{array}{lll}
\renewcommand{\arraystretch}{1.3}
b_\alpha^1 = \{b_{(jkl)}^{i\gamma\delta}\}, &
b_\alpha^2 = \{b_{\alpha kl}^{[\beta\gamma\delta]}\}, &
c_\alpha = \{c_{(ijk)}^{[\alpha\beta\gamma]}\}, \\
b_\beta^1 = \{b_{[jkl]}^{i\gamma\delta}\}, &
b_\beta^2 = \{b_{\alpha kl}^{(\beta\gamma\delta)}\}, &
c_\beta = \{c_{[ijk]}^{(\alpha\beta\gamma)}\}. 
\renewcommand{\arraystretch}{1}
\end{array}
$$

Note that for $p = 2$, we have $b_\alpha^2 = 0$ and 
$c_\alpha = 0$; for $q = 2$, we have $b_\beta^1 = 0$ 
and $c_\beta = 0$; for $p > 2$, we have $c_\alpha = 0$; 
and for $q > 2$, we have  $c_\beta = 0$.

By equations (3.14), the quantities indicated above and the 
subtensors $a_\alpha$ and $a_\beta$ form the following geometric 
objects:
$$
\begin{array}{lll}
\renewcommand{\arraystretch}{1.3}
(a_\alpha, b_\alpha^1),  & (a_\alpha, b_\alpha^2), & 
S_\alpha = (a_\alpha, b_\alpha^1, b_\alpha^2),  \\ 
(a_\beta, b_\beta^1),  & (a_\beta, b_\beta^2), & 
S_\beta = (a_\beta, b_\beta^1, b_\beta^2),  
\renewcommand{\arraystretch}{1}
\end{array}
$$
which are  subobjects of the second structural object and 
the complete structural object of the almost Grassmann 
structure.

Now we consider the cases $p = 2$ and $q = 2$. For definiteness, 
we take the case $p = 2$. As we have already seen, for 
 $p = 2$, the tensor $a_\alpha$ as well as the quantities 
$b_\alpha^2$ and $c_\alpha$ vanish: $a_\alpha = b_\alpha^2 = c_\alpha = 0$, and the object $b_\alpha^1$ becomes a tensor. 
Thus, the vanishing of this tensor is necessary and 
sufficient for the almost Grassmann structure $AG (1, q + 1)$ 
to be $\alpha$-semiintegrable. 

 Hence we have proved the following result.

\begin{theorem} 
\begin{description}
\item[(1)] If $p = 2$, then the structure subobject $S_\alpha$ 
consists only of the tensor $b_\alpha^1$, and the vanishing of 
this tensor is necessary and sufficient for the almost Grassmann 
structure $AG (1, q + 1)$ to be $\alpha$-semiintegrable. 

\item[(2)] If $q = 2$, then the structure subobject $S_\beta$ 
consists only of the tensor $b_\beta^2$, and the vanishing 
of this tensor is necessary and sufficient for the almost 
Grassmann structure $AG (p-1, p + 1)$ $($which is equivalent to 
the structure $AG (1, p + 1))$ to be $\beta$-semiintegrable. 

\item[(3)] If $p = q = 2$, then the complete structural 
object $S$ consists only of the tensors $b_\alpha^1$ and 
$b_\beta^2$, and the vanishing of one of these tensors is 
necessary and sufficient for the almost Grassmann structure 
$AG (1, 3)$ to be $\alpha$- and $\beta$-semiintegrable, 
respectively. 
\rule{3mm}{3mm}
\end{description}
\end{theorem}

We will make two more remarks:

\begin{description}
\item[1)] The tensors $b_\alpha^1$ and $b_\beta^2$ are defined in 
a third-order differential neighborhood of the  almost Grassmann 
structure.

\item[2)]  As was indicated earlier, for $p = q = 2$, the almost 
Grassmann structure $AG (1, 3)$ is equivalent to the conformal 
$CO (2, 2)$-structure. Thus, by results of subsection {\bf 1.5}, 
we have the following decomposition of its 
complete structural object: $S = b_\alpha^1 \dot{+} b_\beta^2$. 
This matches the splitting of the tensor of conformal curvature 
of the $CO (2, 2)$-structure: $b = b_\alpha \dot{+} b_\beta$. 
\end {description}

If $p = 2$ or $q = 2$, the conditions for an almost Grassmann 
structure to be semiflat or flat are connected with a 
differential neighborhood of third order and cannot be expressed 
in terms of the torsion tensor $a$. The reason for this is that 
if $q = 2$, then, as we noted earlier,  the tensor $a_\beta$ 
automatically vanishes, and if $p = 2$, then the same is true for 
the tensor $a_\alpha$. 

Finally, if $p = q = 2$, then for the almost Grassmann structure 
$AG (1, 3)$ we have $a = 0$, and the conditions for an 
 $AG (1, 3)$-structure to be semiflat or flat 
are connected with a differential neighborhood of third order of 
the manifold.

In fact, if $p = q = 2$, then  the almost Grassmann structure 
$AG (1, 3)$ is equivalent to the $CO (2, 2)$-structure. As equations (1.11) show, the torsion tensor of the latter structure 
vanishes, and the conditions for a $CO (2, 2)$-structure to be 
semiintegrable or integrable are expressed in terms of the tensor 
of conformal curvature which is defined in  a differential 
neighborhood of third order of the manifold.

{\bf 2.} The following  table  presents the similarities and 
differences between the conformal structures $CO (p, q)$ 
and almost Grassmann structures $AG (m, n)$:
\vspace*{2mm}

\begin{tabular}{||l|l|l|l||} \hline \hline
 \#& {\em Property} & $CO (p, q)$ & $AG (m, n)$ \\ \hline \hline 
1. & $\dim M$ & $n = p + q$  & $n = p \cdot q, $ \\ &&&$p = m+1, \; q = n-m$  \\ \hline
2. & {\em Invariant const-}  & 2nd order cone 
  & Segre cone  \\ &{\em  ruction in $T_x (M)$}&
$C_x (p, q)$&$SC_x (p, q)$ \\  \hline
3. & {\em Order of} & $s = 1$  & $s = 1$ \\ 
&{\em $G$-structure}&& \\ 
\hline
4. & {\em Structure group} & $G \cong {\bf SO} (p, q) 
\times {\bf H}$  & 
$G \cong {\bf SL} (p) \times {\bf SL} (q) \times {\bf H}$ \\ &&& \\ \hline
5. & {\em Prolonged } & $G' \cong G \htimes {\bf T} (p + q)$  & $G' \cong G 
\htimes {\bf T} (p \cdot q)$  \\ 
&{\em structure group}&& 
\\ \hline
6. & {\em Type of} & $t = 2$  & $t = 2$ \\ 
&{\em $G$-structure}&& \\ \hline
7. & {\em Existence of} & torsion-free &  torsion exists\\ &{\em torsion}&& \\ \hline
8. & {\em Complete } & $\{b, c\}$  & $\{a, b, c\}$ \\ 
&{\em structure  object}&& \\  \hline 
9. & {\em Local space} & $(C^n_q)_x$  &  $G (m, n)_x$\\ &&&\\  \hline
10. & {\em Locally flat} & $ C^n_q$  & $G (m, n)$ \\ & 
{\em structure}&& \\  \hline
11. & {\em Existence of} & $p=q=2$:      & for all $p$ and $q$: 
   \\ &{\em isotropic bundles}& $E_\alpha (M, {\bf SL} (2))$ 
and  & $E_\alpha (M, {\bf SL} (p))$ and \\
&& $E_\beta (M, {\bf SL} (2))$ & $ E_\beta (M, {\bf SL} (q))$ 
\\ \hline
12. & $p+q = p\cdot q 
\Longrightarrow$  & $CO (2, 2) $  & $AG (1, 3)$ \\ 
&$ p=q=2$&& \\ 
\hline \hline
\end{tabular}


\begin{thebibliography}{AHS 78}

\bibitem[A 69]{A:[A 69]} M.~A. Akivis,
   Three-webs of multidimensional surfaces,  
 {\em Trudy  Geom. Sem. Inst. Nauchn. Inform., Akad. Nauk SSSR} 
       {\bf 2} (1969) 7--31 (in Russian).  

\bibitem[A 82a]{A:[A 82a]} M.~A. Akivis,  
 Webs and almost-Grassmann structures,  
{\em Sibirsk. Mat. Zh.} {\bf 23} (1982) (6) 6--15 (Russian).  
English transl. in {\em Siberian Math. J.} 
{\bf 23} (1982) (6) 763--770. 


\bibitem[A 82b]{A:[A 82b]} M.~A. Akivis, 
 On the differential geometry of a Grassmann manifold.  
 {\em Tensor (N.S.)} {\bf 38} (1982) 273--282 
(in Russian). 


\bibitem[A 83]{A:[A 83]} M.~A. Akivis, 
  Completely isotropic submanifolds of a 
four-dimensional  pseudoconformal structure,   
{\em Izv. Vyssh. Uchebn.  Zaved. Mat.} {\bf 1983} (1) (248) 
3--11 (in Russian). 
English transl. in {\em Soviet Math. (Iz. VUZ)} {\bf 27} (1983) 
(1) 1--11.  

\bibitem[A 85]{A:[A 85]} M.~A. Akivis, 
 On the theory of conformal structures,  
 {\em Geom. Sb. Vyp. {\bf 26}} 
(Tomsk. Univ., Tomsk, 1985) 44--52 (in Russian) . 



\bibitem[AG 96]{AG:[AG 96]} M.~A.  Akivis, and V.~V. Goldberg,
 {\em Conformal differential geometry and its generalizations} 
(Wiley-Interscience Publication, New York,   1996).

\bibitem[AK 93]{AK:[AK 93]} M.~A.  Akivis, and V.~V. Konnov,
  Local aspects in conformal structure theory, 
 {\em Uspekhi Mat. Nauk} {\bf 48} (1993),  (1) 3--40 (in Russian).  
English transl. in {\em Russian Math. Surveys} 
{\bf 48} (1993) (1) 1--35.

\bibitem[AS 92]{AS:[AS 92]} M.~A. Akivis, and A.~M. Shelekhov,
{\em  Geometry and algebra of  multidimensional three-webs}  
(Kluwer Academic Publishers, Dordrecht, 1992). 
 

\bibitem[AHS 78]{AHS:[AHS 78]} M.~F. Atiyah,  N.~L. Hitchin,  and 
I. Singer, Self-duality in four-dimensional Riemannian 
 geometry, {\em Proc. Roy. Soc. London Ser. A} 
{\bf 362} (1978) 425--461.  

\bibitem[BE 91]{BE:[BE 91]} T.~N. Bailey,  and M.~G. Eastwood,  
Complex paraconformal manifolds: their differential 
geometry and twistor theory, 
{\em Forum  Math.} {\bf 3} (1991) (1) 61--103. 

\bibitem[B 91]{B:[B 91]} R.~J. Baston,  
Almost Hermitian symmetric manifolds. I. Local twistor 
theory, {\em Duke Math. J.} {\bf 63} (1991) (1) 81--112. 


\bibitem[C 08]{C:[C 08]} \'{E}. Cartan,  
  Les subgroupes des groupes continus de transformations, 
{\em Ann. Sci. \'{E}cole Norm.} (3) {\bf 25} (1908) 57--194; 
{\em \OE uvres compl\`{e}tes. Partie II,  Alg\'{e}bre. 
Formes  Diff\'{e}rentielles, Syst\`{e}mes Diff\'{e}rentielles},  
Vols. 1--2 (Gauthier-Villars, Paris, 1953), 719--856.



\bibitem[C 22a]{C:[C 22a]} \'{E}. Cartan,  
 Sur les \'{e}quations de la gravitation d' Einstein, 
{\em J. Math. Pures Appl.}  {\bf 1}  (1922) 141--203; 
{\em \OE uvres compl\`{e}tes}. Partie III, {\em Divers, 
G\'{e}om\'{e}trie Diff\'{e}rentielle},  
Vols. 1--2, Gauthier-Villars, Paris, 1955, 549--611.

 


\bibitem[C 22b]{C:[C 22b]} \'{E}. Cartan,   
 Sur les espaces conformes g\'{e}n\'{e}ralis\'{e}s 
    et l'Univers optique, {\em C. R. Acad. Sci. Paris} {\bf 174} 
(1922)   857--859; {\em \OE uvres compl\`{e}tes}.
Partie III, {\em Divers, G\'{e}om\'{e}trie diff\'{e}rentielle},  
Vols. 1--2, Gauthier-Villars, Paris, 1955, 622--624.


\bibitem[C 23]{C:[C 23]} \'{E}. Cartan, 
 Les espaces \`{a} connexion conforme, 
{\em  Ann.  Soc. Polon. Math.} {\bf 2} (1923) 171--221;  
{\em \OE uvres compl\`{e}tes. Partie III,  Divers, 
G\'{e}om\'{e}trie Diff\'{e}rentielle},  
Vols. 1--2 (Gauthier-Villars, Paris, 1955) 747--797.

\bibitem[Dh 94]{Dh:[Dh 94]} P.~F. Dhooghe,   
 Grassmannian structures on manifolds, 
{\em Bull. Belg. Math. Soc. Simon Stevin} {\bf 1} (1994) (1) 
597--622. 


\bibitem[E 26]{E:[E 26]}  L.~P. Eisenhart,  
{\em Riemannian geometry} (Princeton Univ. Press, 
Princeton, N.J., 1926;  6th printing,  1966). 


\bibitem[Ga 89]{Ga:[Ga 89]}  R. Gardner, 
{\em The Method of Equivalence and Its Applications} 
(SIAM, Philadelphia, PA, 1989).


\bibitem[G 73]{G:[G 73]}  V.~V. Goldberg,  
 $(n+1)$-webs of multidimensional surfaces,  {\em Dokl.
 Akad. Nauk SSSR} {\bf 210} (1973) (4) 756--759 (in Russian). 
English transl. in {\em Soviet Math. Dokl.} 
{\bf 14} (1973) (3) 795--799.  


\bibitem[G 75]{G:[G 75]} V.~V.  Goldberg,  
 The almost Grassmann manifold that is  connected with an 
 $(n+1)$-web of multidimensional surfaces,  {\em Izv. Vyssh. 
 Uchebn. Zaved. Mat.} {\bf 1975} (8) (159) 29--38 (in Russian).  
 English  transl. in {\em Soviet Math. (Iz. VUZ)} 
{\bf 19} (1975) (8)  23--31.  


\bibitem[G 88]{G:[G 88]} V.~V.  Goldberg, 
{\em Theory of Multicodimensional $(n+1)$-Webs} (Kluwer Academic
 Publishers, Dordrecht, 1988). 


\bibitem[Go 87]{Go:[Go 87]}  A.~B. Goncharov,  
Generalized conformal structures on manifolds, 
{\em Selecta Math. Soviet.}  {\bf 6} (1987)  306--340. 


\bibitem[H 66]{H:[H 66]} Th. Hangan, 
  G\'{e}om\'{e}trie diff\'{e}rentielle grassmannienne, 
{\em Rev. Roumaine Math. Pures Appl.}
 {\bf 11} (1966) (5) 519--531. 


\bibitem[H 68]{H:[H 68]}           
Th. Hangan, Th.,  Tensor-product tangent bundles, 
{\em Arch.  Math. (Basel)} {\bf 19} (1968) (4) 436--440.


\bibitem[H 80]{H:[H 80]}           
Th. Hangan,  Sur l'int\'{e}grabilit\'{e} des structures 
       tangentes produits tensoriels r\'{e}els, 
{\em Ann. Mat. Pura Appl.} (4) {\bf 126} (1980) 149--185. 

\bibitem[I 72]{I:[I 72]}           
T. Ishihara, On tensorproduct structures and Grassmannian 
structures, {\em J. Math. Tokushima Univ.} {\bf 1972} (4) 
1--17. 

\bibitem[M 78]{M:[M 78]}   Yu.~I.      Mikhailov,  
 On the structure of almost Grassmannian manifolds, 
   {\em Izv. Vyssh. Uchebn. Zaved. Mat.} {\bf 1978}, 
(2) 62--72 (in Russian). English  transl. in {\em 
Soviet  Math.    (Iz. VUZ)} {\bf 22} (1978) (2) 54--63. 

\bibitem[R 59]{R:[R 59]}   B.~A. Rosenfeld, 
Quasielliptic spaces 
{\em  Trudy  Moskov. Mat. Obshch.}  {\bf 8} (1959) 
 49--70 (in Russian) .

\bibitem[S 64]{S:[S 64]}  S.  Sternberg,
{\em Lectures on differential geometry} 
(Prentice-Hall, Inc., 
Englewood Cliffs, N.J., 1964; 
2nd edition, Chelsea Publishing Co., New York, N.Y., 1983).  


\bibitem[W 18]{W:[W 18]} H. Weyl, 
Reine Infinitesimalgeometrie, {\em  Math. Z.} {\bf 2} (1918)
 384--411.  

\bibitem[Y 39]{Y:[Y 39]} K. Yano, 
 Sur les equations de Codazzi dans la 
g\'{e}om\'{e}trie conforme des espaces de Riemann, {\em Proc. 
Imp. Acad. Tokyo} {\bf 15} (1939) 340--344.

\end{thebibliography}
\end{document}